\documentclass{article}

\usepackage{arxiv}

\usepackage{graphicx}
\usepackage{amssymb}
\usepackage{mathtools}
\usepackage[lithuanian, english]{babel}
\usepackage[utf8]{inputenc}
\usepackage[L7x]{fontenc}
\usepackage{lmodern}

\DeclareMathOperator{\supp}{supp}
\def \PP{{\bf P}}                 
\def \EE{{\bf E}}                 

\def\H0{\mathrm{H}_0}

\newcommand{\sign}{\text{sign}}

\usepackage[numbers,sort&compress]{natbib}
\bibpunct[, ]{[}{]}{,}{n}{,}{,}
\makeatletter
\def\NAT@def@citea{\def\@citea{\NAT@separator}}
\makeatother
\usepackage[caption=false]{subfig}

\usepackage{enumitem}
\usepackage{hyperref}

\newtheorem{theorem}{Theorem}[section]

\newtheorem{remark}{Remark}

\usepackage{hyperref}
\hypersetup{
    colorlinks=true,
    linkcolor=blue,
    filecolor=magenta,      
    urlcolor=blue,
    citecolor = blue
}

\begin{document}

\title{Multiple outlier detection tests for parametric models}

\author{
	Vilijandas Bagdonavi\v{c}ius \\
	Institute of Applied Mathematics,\\ Vilnius University,  Vilnius, Lithuania\\
	\texttt{vilijandas.bagdonavicius@mif.vu.lt} \\
	\And
	Linas Petkevi\v{c}ius \\
	Institute of Computer Science,\\ Vilnius University,  Vilnius, Lithuania \\
	\texttt{linas.petkevicius@mif.vu.lt}
}


\maketitle

\begin{abstract}
We propose a simple multiple outlier identification method for parametric location-scale and shape-scale 
models when the number of possible outliers is not specified. The method is based on a result giving  asymptotic properties of extreme z-scores. Robust estimators of model parameters are used defining z-scores. An extensive simulation study was done for comparing of  the proposed method with existing methods.  For the normal
family, the method is compared with the well known Davies-Gather, 
Rosner’s, Hawking's and Bolshev's multiple outlier identification methods. The choice of an upper limit for the number of possible outliers in case of Rosner's test 
application is discussed. For other  families, the proposed method  is  compared
with a method generalizing Gather-Davies method. In most situations, the new method has the highest 
outlier identification power in terms of masking and swamping values. We also created R package outliersTests for proposed test.
\end{abstract}
{
\small
  \textbf{\textit{Keywords---}}
Location-scale  models; Outliers identification; Unknown number of outliers; Outlier region, Robust estimators
}

\section{Introduction}

The problem of multiple outliers identification received attention of many
authors. The majority of  outlier identification methods define  rules  for the rejection of the most extreme observations.
The bulk of publications  have been concentrated on the normal distribution (see \cite{bolshev_chauvenets_1975,davies_identification_1993,dixon_analysis_1950,grubbs_sample_1950,rosner_detection_1975,tietjen_grubbs-type_1972},  see surveys in \cite{barnett1974outliers,zerbet_statistical_2002}. For non-normal case,
the most of the literature pertains to the exponential and gamma distributions, see
\cite{chikkagoudar_distributions_1983,kabe_testing_1970,kimber_testing_1988,lalitha_multiple_2012,lewis_recursive_1979,likes_distribution_1967,lin_exact_2009,lin_tests_2014,zerbet_new_2003}.

Constructing outlier identification methods, the most of authors suppose that the number $s$ of  observations suspected to be outliers  is specified. These methods have a serious drawback: 
only two possible conclusions are done:  exactly $s$ observations are admitted as  outliers or it is concluded that outliers are absent. More natural is to consider methods 
which do not specify  the number of
 suspected observations or at least specify the upper limit $s$ for it. Such methods are not very numerous and they concern  normal or exponential samples. These are  \cite{rosner_detection_1975,bolshev_chauvenets_1975,hawkins_identification_1980} methods for normal samples, 
 \cite{kimber_tests_1982,lin_exact_2009,lin_tests_2014}) methods for exponential samples. The only method which does not specify the upper limit $s$  is the  \cite{davies_identification_1993} method for normal samples. 

We give a competitive and simple method for outlier identification in samples from location-scale and shape-scale families of probability distributions. The upper limit $s$ is not specified, as in the the case of Davies-Gather method.   The method is based on a theorem giving  asymptotic properties of extreme z-scores. Robust estimators of model parameters are used defining z-scores.

In Section \ref{sec:outliers}  we present a short  overview of the notion of the outlier region given by  \cite{davies_identification_1993}. In Section \ref{sec:new_method} we give asymptotic properties of extreme z-scores based on equivariant estimators of model parameters, and
 introduce  a new outlier identification method for  parametric models based on  the asymptotic result and robust  estimators.  In section \ref{sec:davies} we consider rather evident generalizations of Davies-Gather tests for normal data to location-scale families. In Section \ref{sec:compare} we give a short overview of known multiple outlier identification methods for normal samples which do not specify an exact number of suspected outliers. 
 In Section \ref{sec:simul} we compare performance of the new and existing methods.

\section{Outliers and outlier regions}
\label{sec:outliers}
Suppose that  data are independent random variables
$X_1,\ldots,X_n$. Denote by $F_{i}(x)$ the c.d.f. of $X_i$.

Let
${\cal F}_0=\{F(x,\theta ), \theta \in \Theta \subset
\mathbf{R}^m\} \label{eq:param_family}
$
be a parametric family of absolutely continuous cumulative distribution functions with continuous unimodal densities $f$  on the support
$\supp(F)$  of the c.d.f. $F$.

\vskip 0.2cm \noindent Suppose that if the data are not contaminated with unusual observations, then  the following  null  hypothesis $H_0$ is true: there exist $\theta\in \Theta$ such that
\begin{equation}
F_{1}(x)=\ldots=F_{n}(x)=F(x,\theta).\label{eq:null hypothesis}
\end{equation}

There are two different definitions of an outlier. In the first case
outlier is an observation which falls  into some outlier region $out(X)$. The outlier region is a set such   that the probability
for at least one observation from a sample to  fall into it is small if the hypothesis $H_0$ is true. In such a case the probability that a specified observation $X_i$  falls into $out(X)$ is very small. If an observation $X_i$ has  distribution different from that  under $H_0$ then this probability may be considerably higher. 

In the second case the value $x_i$ of $X_i$ is an outlier if the probability distribution of $X_i$ is different from that under $H_0$, formally $F_i \neq F(x,\theta)$. In this case outliers are often called {\it contaminants}. 

So in the first case exists a very small probability to have
an outlier under $H_0$. In the second case outliers (contaminants) do not exist under $H_0$ and  outliers (contaminants)  do not necessary fall into the outlier region. Both definitions give approximately the same outliers if the alternative distribution is concentrated in the outlier region. Namely such contaminants can be called outliers in the sense that outliers are anomalous extreme observations. In such a case it is possible to compare outlier and contaminant search methods.

In this paper, we consider  location-scale and shape-scale families. Location-scale families have the form ${\cal F}_{ls}=\{F_0((x-\mu)/\sigma), \mu \in {\bf R},\sigma>0\}$   with the completely specified baseline
c.d.f $F_0$ and  p.d.f. $f_0$. Shape-scale families have the form ${\cal F}_{ls}=\{G_0(((x/\theta)^\nu), \theta, \nu >0\}$   with completely specified baseline
c.d.f $G_0$ and  p.d.f. $g_0$. By logarithmic transformation the shape-scale families are transformed to  location-scale family, so we concentrate on location-scale families. Methods for such families are easily modified to methods for shape-scale families.   

The right-sided $\alpha$-outlier region for a location-scale family is
\begin{equation*}
out_r(\alpha_n,F)=\{x\in {\bf R}:x>\mu+\sigma F_0^{-1}(1-\alpha)\}
\end{equation*}
and the left-sided $\alpha$-outlier region is
\begin{equation*}
out_l(\alpha_n,F)=\{x\in {\bf R}:x<\mu+\sigma F_0^{-1}(\alpha)\}.
\end{equation*}
The two-sided $\alpha$-outlier region has the form
\begin{equation}
out(\alpha,F)=\{x\in {\bf R}/[\mu+\sigma F_0^{-1}(\alpha/2),\mu+\sigma F_0^{-1}(1-\alpha/2)]\}.
\end{equation}

If $f_0$ is symmetric, then the two-sided outlier region is simpler:
\begin{equation*}
out(\alpha,F)=\{x\in {\bf R}:|x-\mu|>\sigma F_0^{-1}(1-\alpha/2)\}.
\end{equation*}
The value of $\alpha$
is chosen depending on the size $n$ of a  sample:
$\alpha=\alpha_n$. The choice is based on assumption that
under $H_0$ for some $\bar\alpha$  close to zero
\begin{eqnarray}
\PP\{\cap_{i=1}^n\{X_i \notin out(\alpha_n,F)\}\}= (\PP\{X_i
\notin out(\alpha_n,F)\})^n =1-\bar \alpha. \label{eq:out_reg_p}
\end{eqnarray}
 The equality
\eqref{eq:out_reg_p} means that under $H_0$ the probability that
{\it none} of $X_i$ falls into $\alpha_n$ - outlier region is
$1-\bar\alpha$. It implies that
\begin{equation} \alpha_n=1-(1-\bar\alpha)^{1/n}.\label{eq:alpha n}\end{equation} The sequence
$\alpha_n$ decreases from $\bar \alpha$ to  $0$ as $n$ goes from
$1$ to $\infty$.

The first definition of an outlier is as follows: for a sample size $n$ a realization $x_i$ of $X_i$ is called
{\it outlier} if $x_i\in out(\alpha_n,F)$;   $x_i$  is called
{\it right  outlier} if $x_i\in out_r(\alpha_n,F)$.

The number of outliers $D_n$ under $H_0$ has the binomial
distribution $B(n,\alpha_n)$ and the expected number of outliers
in the sample under $H_0$ is $\EE D_n=n\alpha_n.$ Note that $\EE D_n\to
-\ln(1- \bar\alpha)\approx \bar\alpha$ as $n\to\infty$. For example, if
 $\bar\alpha=0.05$, then $\ln(1-\bar\alpha)\approx
0.05129$ and for  $n\geq 10$ the expected number of outliers  is
approximately $0.051$, i.e. it practically does not depend on $n$. So  under $H_0$  the expected
number of outliers  $0.051$ is  negligible with respect to the sample size $n$.

\section{New method}
\label{sec:new_method}
\subsection{Preliminary results}

Suppose that a c.d.f. $F\in {\cal F}_{ls}$ 
belongs also to the domain of attraction ${\cal G}_\gamma$, $\gamma\geq  0$ (see \cite{de2007extreme}).

If $F\in {\cal G}_0\cap {\cal F}_{ls}$ , then  there exist normalizing constants $a_n>0 $
and $b_n\in {\bf R}$ such that $\lim_{n\to\infty} F_0^n(a_n x+b_n) = e^{-e^{-x}}.$  Similarly, if 
$F\in {\cal G}_\gamma\cap {\cal F}_{ls}$, $\gamma>0$,  then $\lim_{n\to\infty} F_0^n(a_n x+b_n) =  e^{-(-x)^{-1/\gamma}}, \quad x<0$, $\lim_{n\to\infty} F^n(a_n x+b_n) = 1, x\geq 0.$

One of possible choices of the sequences  $\{b_n\}$ and $\{a_n\}$ is
\begin{equation}
 b_n=F_0^{-1}(1-\frac1n),\quad a_n=1/(nf_0(b_n)). \label{eq:anbn}
\end{equation}
In the particular case of the normal distribution equivalent form $a_n=1/b_n$ can be used. Expressions of $b_n$ and $a_n$ for some most used distributions are given in Table \ref{tab:an_bn}.

\begin{table}
 \caption{Expressions of $b_n$ and $a_n$}\label{tab:an_bn}
 \centering 
\begin{tabular}{|c|c|c|c|}
  \hline
  Distribution& $F_0(x)$ & $b_n$ & $a_n$  \\
  \hline
  Normal&$\Phi(x)$ &$\Phi^{-1}(1-1/n)$ &$1/b_n$ \\
  Type I extreme value&$1-e^{-e^{x}}$ &$\ln\ln n$&$e^{-b_n}$ \\
  Type II extreme value& $e^{-e^{-x}}$  &$\ln(-\ln(1-1/ n))$ & $e^{b_n}/(n-1)$ \\
  Logistic&$ \frac1{1+e^{-x}}$  &$  \ln(n-1)$ &$ {n}/({n-1})$ \\
  Laplace & $\frac12+\frac12\sign(x)(1-e^{-|x|})$ & $ \ln({n}/2)$& $1$ \\
  Cauchy & $\frac12+\frac1{\pi}\arctan (x)$  & $\cot(\frac\pi{n})$ & $\frac{\pi}{n}/\sin^2(\frac\pi{n})$
\\
\hline
\end{tabular}

\end{table}

\begin{cond} Consider a model that satisfies the following conditions:\label{cond:cond}
\begin{enumerate}[label=\alph*)]
\item  $\hat\mu$ and $\hat \sigma$ are consistent estimators of $\mu$ and $\sigma$;\label{cond:a}

\item the limit distribution of ($\sqrt{n}(\hat\mu-\mu), \sqrt{n}(\hat\sigma-\sigma))$ is non-degenerate;\label{cond:b}

\item  
\begin{equation*}
    \lim_{x\to \infty} \frac{xf_0(x)}{\sqrt{1-F_0(x)}}= 0.
\end{equation*}\label{cond:c}
\end{enumerate}
\end{cond}

\hyperref[cond:cond]{Condition A}  \ref{cond:c} is satisfied for many location-scale models including  the normal, type I extreme value, type II extreme value, logistic, Laplace ($F\in {\cal G}_0$), Cauchy
($F\in {\cal G}_1$).

Set $Y_i=(X_i-\mu)/{\sigma}$, $\hat Y_i=(X_i-\hat\mu)/{\hat\sigma}$. The random variables $\hat Y_i$ are called z-scores. Denote by $Y_{(1)}\leq Y_{(2)}\leq \ldots \leq Y_{(n)}$ and $\hat Y_{(1)}\leq \ldots \leq \hat Y_{(n)}$ the respective order statistics
 
 The following theorem is useful for  right outliers detection test construction.
\begin{theorem}\label{theo:teorem1}

If  $F\in {\cal G}_0\cap {\cal F}_{ls}$ and \hyperref[cond:cond]{Conditions A} hold, then for fixed
$s$
$$
((\hat Y_{(n)}-b_n)/{a_n},(\hat Y_{(n-1)}- b_n)/{ a_n}
,...,(\hat Y_{(n-s+1)}- b_n)/{ a_n})\stackrel{d}\rightarrow$$$$
L_0=(-\ln E_1,-\ln (E_1+E_2),...,-\ln (E_1+...+E_s))
$$
as $
n\to\infty$, where $E_1,...,E_s$ are i.i.d. standard
exponential random variables.

If $F\in {\cal G}_\gamma\cap {\cal F}_{ls}$, $\gamma> 0$ and \hyperref[cond:cond]{Conditions A}  hold, then the limit random vector is
$$
L_\gamma=( E_1^{-1}-1,(E_1+E_2)^{-1}-1,...,(E_1+...+E_s)^{-1}-1).
$$

\end{theorem}

\vskip 0.2cm
\begin{remark}

Note that  $2(E_1+...+E_i)\sim \chi^2(2i)$.
  It implies that if $F\in {\cal G}_0\cap {\cal F}_{ls}$, then for fixed $i$, $i=1,...,s$,
\begin{equation}
\PP\{(\hat Y_{(n-i+1)}-b_n)/a_n\leq x\}\to 
1-F_{\chi^2_{2i}}(2e^{-x})\quad \mbox{as}\;\;n\to\infty.
\end{equation}
Similarly, if  $F\in {\cal G}_\gamma\cap {\cal F}_{ls}$, $\gamma>0$,  then for fixed $i$, $i=1,...,s$,
\begin{equation}
\PP\{(\hat Y_{(n-i+1)}-b_n)/a_n\leq x\}\to
1-F_{\chi^2_{2i}}(\frac2{1+x
})\quad \mbox{as}\;\;n\to\infty.
\end{equation}

\end{remark}

The following theorem is useful for construction of outlier detection tests
in two-sided case when  $f_0$ is symmetric.  For any sequence $\zeta_1,...,\zeta_n$ denote by $|\zeta|_{(1)}\leq ... \leq |\zeta|_{(n)}$ the ordered absolute values $|\zeta_1|,..., |\zeta_n|$. 
\begin{theorem}\label{theo:teorem2}

Suppose that the function $f_0$ is symmetric. If  $F\in {\cal G}_\gamma\cap {\cal F}_{ls}$, $\gamma\geq 0$ and \hyperref[cond:cond]{Conditions A} hold, then for fixed
$s$
$$
((|\hat Y|_{(n)}-b_{2n})/{a_{2n}},(|\hat Y|_{(n-1)}-b_{2n})/{a_{2n}}
,...,(|\hat Y|_{(n-s+1)}-b_{2n})/{a_{2n}})\stackrel{d}\rightarrow L_\gamma
$$
as $
n\to\infty$. 
\end{theorem}
\vskip 0.2cm


\begin{remark}\label{theo:remark2}
\hyperref[theo:teorem2]{Theorem 3.2} implies that if $F\in {\cal G}_0\cap {\cal F}_{ls}$, $n\to\infty$, then for fixed $i$, $i=1,...,s$,
\begin{equation}
\PP\{(|\hat Y|_{(n-i+1)}-b_{2n})/a_{2n}\leq x\}\to 
1-F_{\chi^2_{2i}}(2e^{-x}),
\end{equation}
and if $F\in {\cal G}_\gamma\cap {\cal F}_{ls}$, $\gamma>0$,  then
\begin{equation}
\PP\{(|\hat Y|_{(n-i+1)}-b_{2n})/a_{2n}\leq x\}\to 1-F_{\chi^2_{2i}}(2/(1+x)).
\end{equation}
\end{remark}
Suppose now that the function $f_0$ is not symmetric. Set  $Y^*_i=-(X_i-\mu)/\sigma$. The c.d.f. and p.d.f. of $Y^*_i$ are $1-F_0(-x)$
 and $f_0(-x)$, respectively. Set \begin{equation}\label{eq:normal_seq_star}
  b^*_{n}=-F_0^{-1}(1/n), \quad a^*_{n}=1/({nf_0(-b^*_{n})}).
\end{equation}
 For example, if type I extreme value distribution is considered, then 
 $$
  b_{n}=\ln\ln n, \quad a_{n}=\frac1{\ln n},\quad   b^*_{n}=-\ln(-\ln(1-\frac1n)), \quad a^*_{n}=-\frac1{(n-1)\ln(1-\frac1n)}.
 $$
 For the type II extreme value distribution $a_n,b_n,a_n^*,b_n^*$ have the same expressions as $a^*_n,b^*_n,a_n,b_n$ 
 for the Type I extreme value distribution, respectively.

 \begin{remark}\label{theo:remark3}
 Similarly as in \hyperref[theo:teorem1]{Theorem 3.1} we have that if $s$ is fixed and $F\in {\cal G}_0\cap {\cal F}_{ls}$, then 
  for fixed $i$, $i=1,...,s$,
\begin{equation}
\PP\{( Y_{(i)}+b^*_{n})/(-a^*_{n})\leq x\}=\PP\{(\hat Y^*_{(n-i+1)}-b^*_{n})/a^*_{n}\leq x\}\to
1-F_{\chi^2_{2i}}(2e^{-x}),
\end{equation}
and if  $F\in {\cal G}_\gamma\cap {\cal F}_{ls}$, $\gamma> 0$, then 
  for fixed $i$, $i=1,...,s$,
\begin{equation}\label{normal_const_left}
\PP\{( Y_{(i)}+b^*_{n})/(-a^*_{n})\leq x\}=\PP\{(\hat Y^*_{(n-i+1)}-b^*_{n})/a^*_{n}\leq x\}\to
1-F_{\chi^2_{2i}}(2/(1+x)).
\end{equation}
\end{remark}
\subsection{Robust estimators for location-shape distributions}
The choice of the estimators $\hat\mu$ and $\hat \sigma$ is important when outlier detection problem is considered.
The ML estimators from the complete sample  are not stable when
 outliers exist.

In the case of location-scale families highly efficient robust estimators of the location and scale parameters $\mu$ and  $\sigma$ are (see \cite{rousseeuw_alternatives_1993})
\begin{equation}
\hat\mu=MED-\hat\sigma F_0^{-1}(0.5),\quad \hat\sigma=Q_n=d \,W_{([0.25 n(n-1)/2])},\label{ref:robust estimators}
\end{equation}
where $MED$ is the empirical median, $W_{ij}=|X_i-X_j|,\quad 1\leq i<j\leq n$ are $C_n^2=n(n-1)/2$ absolute values of the differences $X_i-X_j$ and $W_{(l)}$ is
the $l$th order statistic from $W_{ij}$.

The constant $d$ has the form
$d=1/K_0^{-1}(5/8)$,
where $K^{-1}_0(x)$ is the inverse of the c.d.f. of $Y_1-Y_2$, $Y_i=(X_i-\mu)/\sigma\sim F_0(x)$.

Expressions of $K^{-1}_0(x)$ and values $d$ for some well-known  location-scale families are given in Table \ref{tab:cn_dn}.
\vskip 0.2cm
\begin{table}
 \caption{Values of $d$  for various probability distributions}\label{tab:cn_dn}
\centering
\begin{tabular}{|c|c|c|c|}
  \hline
  Distribution & $K_0(x)$ &d \\
  \hline
  Normal &$\Phi(x/\sqrt{2})$&2.2219 \quad \\
  Type I extr.val. &$1/(1+e^{-x})$&1.9576   \\
  Type II extr.val.  \quad & $1/(1+e^{-x})$&1.9576  \\
  Logistic   &$ 1-\frac{(x-1)e^x+1}{(e^x-1)^2}$ &1.3079 \\
  Laplace  & $1-\frac12 (1+\frac{x}2)e^{-x}$&1.9306\\
  Cauchy   & $\frac12+\frac1{\pi}\arctan (x/2)$\quad &1.2071 \\
  \hline
\end{tabular}
\end{table}

The above considered estimators are equivariant under $H_0$, i.e.
 for any $e\in {\bf R}, f>0$, the following
equalities hold:
$$
\hat\mu((X_1-e)/f,\ldots,(X_n-e)/f)=(\hat\mu(X_1,\ldots,X_n)-e)/f,
$$$$
\hat\sigma((X_1-e)/f,\ldots,(X_n-e)/f)=\hat\sigma(X_1,\ldots,X_n)/f.
$$
Equivariant estimators have the following property: the distribution of $(\hat\mu-\mu)/\sigma$, $\hat\sigma/\sigma$ and $(\hat\mu-\mu)/\hat\sigma$ is parameter-free.

\subsection{Right outliers identification method for location-scale families}

Suppose that $F\in {\cal G}_\gamma\cap {\cal F}_{ls}$, $\gamma\geq 0$.
 Let $a_n,b_n$ be defined by
 \eqref{eq:anbn}.
Set
$$
 U^+_{(n-i+1)}(n)=1-F_{\chi^2_{2i}}(2e^{-(\hat Y_{(n-i+1)}- b_n)/ a_n}),\;\gamma=0,$$
 $$ U^+_{(n-i+1)}(n)=1-F_{\chi^2_{2i}}(2/(1+(\hat Y_{(n-i+1)}- b_n)/ a_n),\;\gamma>0, $$
\begin{equation} 
U^+(n,s)=\max_{1\leq i\leq s}U^+_{(n-i+1)}(n).
\end{equation}

\begin{theorem}
The distribution of the statistic $U^+(n,s)$ is parameter-free for any fixed $n$
\end{theorem}

\vskip 0.2cm

Denote by $u^+_\alpha(n,s)$ the $\alpha$ critical value of the
statistic $U^+(n,s)$.  Note that it is exact, not asymptotic $\alpha$
critical value: $\PP\{U^+(n,s)\geq u^+_\alpha(n,s)\}=\alpha$ under
$H_0$.

\hyperref[theo:teorem1]{Theorem 3.1} implies that 
 the limit distribution (as $n\to\infty$) of  the random variable $U^+(n,s)$ coincides with  the distribution of the random variable
$V^+(s)=\max_{1\leq i\leq s}V^+_i,$ where $ V^+_i=1-F_{\chi^2_{2i}}(2(E_1+...+E_i))$, $E_1,...,E_s$ are i.i.d. standard exponential  random variables. The random variables $V^+_1,...,V^+_s$ are dependent identically distributed and the distribution of each $V^+_i$ is uniform: $V^+_i\sim U(0,1)$.

Denote by  $v^+_{\alpha}(s)$ the $\alpha$ critical values of the random variable $V^+(s)$. They are easily found by simulation  many times generating $s$ i.i.d. standard exponential random variables and computing values of the random variables $V^+(s)$.

Our simulations showed that the  below proposed outlier identification methods based on exact and approximate critical values of the statistic $U^+(n,s)$  give practically the same results, so
 for samples of size $n\geq 20$ we recommend to  approximate the $\alpha$-critical level of the statistic  $U^+(n,s)$ by 
 the critical values $v^+_\alpha(s)$ which depend only on $s$. We shall see that for the purpose of outlier identification only the critical values $v^+_\alpha(5)$ are needed.  We found that the critical values $v^+_{\alpha}(5)$ are: $v^+_{0.1}(5)=0.9677$,  $v^+_{0.05}(5)=0.9853$,  $v^+_{0.01}(5)=0.9975$.

Our simulations showed that the performances of the below proposed outlier identification method based on exact and approximate critical values of the statistic $U^+(n,5)$  are similar for samples of size $n\geq 20$.

We write shortly BP-method for the below considered method.
\vskip 0.2cm
{\it BP  method for right outliers.} Begin outlier search using observations corresponding to the largest  values of $\hat Y_i$. We recommend begin with five  largest. So take $s=5$ and
compute  the values of the statistics $$U^+(n,5)=\max_{1\leq i\leq 5}U^{+}_{(n-i+1)}(n).$$

 If $U^+(n,5) \leq v^+_\alpha(5)$,  then it is concluded that outliers  are
absent and no further investigation is done. Under $H_0$ the probability of such event is approximately $1-\alpha$.

 If $U^+(n,5) > v^+_\alpha(5)$,  then it is concluded that outliers exist.
 
 Note that (see the classification scheme below) that if $U^+(n,5) > v^+_\alpha(5)$,  then minimum one observation is declared as an outlier. So the probability to declare  absence of outliers does not depend on the following classification scheme.

If it is concluded that outliers exist then search of outliers is done using the following steps.

{\bf Step 1.} Set $d_1=\max\{i\in\{1,...,5\}:U^+_{(n-i+1)}(n)> v^+_\alpha(5)\}$. Note that the maximum $d_1>0$ exists because $U^+(n,5) > v^+_\alpha(5)$.

  If $d_1< 5$, then classification is finished at this step:  $d_1$ observations  are declared as right outliers because if 
the value of  $X_{(n-d_1)}$ is declared as an outlier, then it is natural  to declare values of $X_{(n)},...,X_{(n-d_1+1)}$  as  outliers, too.   

 If $d_1=5$,  then it is possible that the number of outliers is higher than $5$.
 Then the observation corresponding to $i=1$ (i.e corresponding to $X_{(n)}$) is declared as an outlier and we proceed to the step 2.

 {\bf Step 2.} The above written procedure is repeated
 taking $U^+(n-1,5)=\max_{1\leq i\leq 5}U^+_{(n-i)}(n-1)$ instead of $U^+(n,5)$; here
$$
U^+_{(n-i)}(n-1)=1-F_{\chi^2_{2i}}(2e^{-(\hat Y_{(n-i)}- b_{n-1})/ a_{n-1}}),\quad i=1,...,5,
$$
 Set $d_2=\max\{i\in\{1,...,5\}:U^+_{(n-i)}(n-1)> v^+_\alpha(5)\}$.  If  $d_2<5$, the classification is finished  and $d_2 + 1$ observations are declared as  outliers.

 If $d_2=5$,  then it is possible that the number of  outliers is higher than $6$.
 Then the observation corresponding to the largest  $\hat Y_{(n-1)}$ is declared as an outlier, in total $2$ observations (i.e.  the observations corresponding to $i=1,2$ (i.e corresponding to $X_{(n)}$ and $X_{(n-1)}$) are declared as  outliers and we proceed to the step 3. And so on.
 Classification  finishes at the $l$th step when $d_l<5$. So we declare $(l-1)$  outliers in the previous steps and $d_l$  outliers in the last one. The total number of observations declared as  outliers is $l-1+d_l$.
 These observations are values of $X_{(n)},...,X_{(n-d_l-l+2)}$.

\subsection{Left outliers identification method for location-scale families}
Let $a_n^*,b_n^*$ be the normalizing constants defined by \eqref{eq:normal_seq_star}. If $F\in {\cal G}_0\cap {\cal F}_{ls}$,  $i=1,...,s$, then  set
$$
  U^-_{(i)}(n)=1-F_{\chi^2_{2i}}(2e^{(\hat Y_{(i)}+ b^*_n)/ a^*_n}),\quad 
 U^-(n,s)=\max_{1\leq i\leq s}U^-_{(i)}(n).
$$
If $F\in {\cal G}_\gamma\cap {\cal F}_{ls}$, $\gamma> 0$, then replace $e^{(\hat Y_{(i)}+ b^*_n)/ a^*_n}$ by $1/{(1+(\hat Y_{(i)}+ b^*_n)/ a^*_n})$.
Denote by $u^{-}_\alpha(n,s)$   the $\alpha$ critical value of the
statistic $U^-(n,s)$.

\hyperref[theo:teorem1]{Theorem 3.1}  and \hyperref[theo:remark3]{Remark 3}  imply that 
 the limit distribution (as $n\to\infty$) of  the random variable $U^-(n,s)$ coincides with  the distribution of the random variable
$V^+(s)$. So the critical values $u^{-}_\alpha(n,s)$ are approximated by the critical values  $v^-_\alpha(s)=v^+_{\alpha}(s)$.

The left outliers search method coincides with the right outliers search method replacing $+$ to $-$ in all formulas.  
\subsection{Outlier detection tests for location-scale families: two-sided alternative, symmetric distributions}

Let $a_n,b_n$ be defined by \eqref{eq:anbn}.
If $F\in {\cal G}_0\cap {\cal F}_{ls}$,  $i=1,...,s$, then  set
 $$U_{(n-i+1)}(n)=1-F_{\chi^2_{2i}}(2e^{-(|\hat Y|_{(n-i+1)}- b_{2n})/ a_{2n}}),
 \quad U(n,s)=\max_{1\leq i\leq s}U_{(n-i+1)}(n)
.$$
If $F\in {\cal G}_\gamma\cap {\cal F}_{ls}$, $\gamma> 0$, then replace $e^{(\hat Y_{(i)}+ b^*_n)/ a^*_n}$ by $1/{(1+(\hat Y_{(i)}+ b^*_n)/ a^*_n})$. Denote by $u_\alpha(n,s)$ the $\alpha$ critical value of the
statistic $U(n,s)$.

\hyperref[theo:teorem1]{Theorem 3.1} and \hyperref[theo:remark2]{Remark 2} imply that 
 the limit distribution (as $n\to\infty$) of  the random variable $U(n,s)$ coincides with  the distribution of the random variable
$V^+(s)$. So the critical values $u_\alpha(n,s)$ are approximated by the critical values  $v_\alpha(s)=v^+_{\alpha}(s)$.

The outliers search method coincides with the right outliers search method skipping upper index $+$  in all formulas.

 \subsection{Outlier detection tests for location-scale families: two-sided alternative, non-symmetric distributions}

Suppose now that the function $f_0$ is not symmetric. 
  Let $a_n,b_n,a^*_n,b^*_n$ be defined by \eqref{eq:normal_seq_star}.

Begin outlier search using observations corresponding to the largest and the smallest values of $\hat Y_i$. 
We recommend begin with five smallest and 
five largest. So 
compute  the values of the statistics $U^-(n,5)$ and 
$U^+(n,5)$. 
 If $U^-(n,5) \leq v_{\alpha/2}(5)$ and $U^+(n,5) \leq v_{\alpha/2}(5)$,  then it is concluded that outliers  are
absent and no further investigation is done.

 If $U^-(n,5) > v_{\alpha/2}(5)$ or $U^-(n,5) > v_{\alpha/2}(5)$,  then it is concluded that outliers exist. 
 If $U^-(n,5) > v_{\alpha/2}(5)$, then left outliers are searched as in Section 3.3. If $U^+(n,5) > v_{\alpha/2}(5)$, then right outliers are searched as in Section  3.2. The only difference is that $\alpha$ is replaced by $\alpha/2$ in all formulas.

\subsection{Outlier identification method for shape-scale families}
If shape-scale families of the form $\{F(t;\theta,\nu)=G_0((t/\theta)^\nu),\;\theta,\nu> 0\}$ with
specified $G_0$ are considered then the above given tests for location-scale families could be used because
if $X_1,...,X_n$ is a sample from shape scale family then $Z_1,...,Z_n$, $Z_i=\ln X_i$, is a sample from location-scale
family $\{F_0((x-\mu)/\sigma,\mu\in {\bf R}, \sigma>0\})$ with $\mu=\ln \theta$, $\sigma=1/\nu$, $F_0(x)=G_0(e^x)$.

\subsection{Illustrative example}
To illustrate simplicity of the  BP-method, let us consider an illustrative example of its application (sample size  $n=20$, 
$r=7$ outliers).  The sample of size $n=20$ from standard normal distribution  was generated. The 1st-3rd and 17th-20th observations were  replaced by outliers.  The observations $x_i$, the absolute values $|\hat{Y}_i|$ of the z-scores $\hat{Y}_i$, and the ranks $(i)$ of  $|\hat{Y}_i|$ are presented in Table \ref{table:example1}.

\begin{table}[t]
\caption{Illustrative sample ($n=20, r=7$)}\label{table:example1}
\centering 
\begin{tabular}{|rrrrr|rrrr|}
  \hline
 & i & $x_i$ & $|\hat{Y}_i|$ & (i) & i & $x_i$ & $|\hat{Y}_i|$ & (i) \\ 
  \hline
 & 1 & 6.10 & 3.18 & 16 & 11 & -0.69 & 0.28 & 9 \\ 
 & 2 & 10 & 5.17 & 18 & 12 & -0 & 0.07 & 5 \\ 
 & 3 & 6.20 & 3.23 & 17 & 13 & 0.05 & 0.10 & 6 \\ 
 & 4 & -0.08 & 0.03 & 2 & 14 & -0.20 & 0.03 & 1 \\ 
 & 5 & 0.63 & 0.39 & 11 & 15 & -0.25 & 0.06 & 4 \\ 
 & 6 & -0.54 & 0.21 & 7 & 16 & -0.64 & 0.25 & 8 \\ 
 & 7 & 1.37 & 0.77 & 13 & 17 & -6.30 & 3.14 & 15 \\ 
 & 8 & 0.46 & 0.30 & 10 & 18 & -5.50 & 2.73 & 14 \\ 
 & 9 & -0.22 & 0.04 & 3 & 19 & -12.10 & 6.10 & 19 \\ 
 & 10 & 0.94 & 0.55 & 12 & 20 & -20 & 10.13 & 20 \\ 
   \hline
\end{tabular}
\end{table}

\begin{table}[t]
\caption{Illustrative example of BP test observations classification.}\label{table:example11}
\centering 
\begin{tabular}{|rrrrrr|r|}
  \hline
 & $U_{(20)}(20)$ & $U_{(19)}(20)$ & $U_{(18)}(20)$ & $U_{(17)}(20)$ & $U_{(16)}(20)$ & $U(20, 5)$\\ 
 & \bf  1.000000 & 1.000000 & 1.000000 & 0.999998 & 1.000000  & 1.000000\\ \hline
  & $U_{(19)}(19)$ & $U_{(18)}(19)$ & $U_{(17)}(19)$ & $U_{(16)}(19)$ & $U_{(15)}(19)$ & $U(19, 5)$\\
  & \bf 0.999685 & 0.999998 & 0.999916 & 0.999998 & 1.000000 & 1.000000\\  \hline 
    & $U_{(18)}(18)$ & $U_{(17)}(18)$ & $U_{(16)}(18)$ & $U_{(15)}(18)$ & $U_{(14)}(18)$ & $U(18, 5)$\\
  & \bf 0.998046 & 0.996970 & 0.999893 & 0.999997 & 0.999997 & 0.999997\\ \hline
    & $U_{(17)}(17)$ & $U_{(16)}(17)$ & $U_{(15)}(17)$ & $U_{(14)}(17)$ & $U_{(13)}(17)$ & $U(17, 5)$\\
  &\bf 0.924219 &\bf  0.996446 &\bf  0.999871 &\bf  0.999940 & 0.084290 & 0.999940\\ \hline
   \hline
\end{tabular}
\end{table}

In Table \ref{table:example11} we present steps of the classification procedure by the BP method. First, we compute (see line 1 of Table \ref{table:example11})  value of the statistic\footnote{In fact the value of the statistic is 0.999999999.  It is rounded to 1.} $U(20, 5)= \max_{1\leq i\leq 5}U_{(20-i+1)}(20)= 1$. Since $U(20, 5)=1 > 0.9853 = v_{0.05}(5)$,  we reject the null hypothesis, conclude that outliers exist and begin the search of outliers.

{\it Step 1}.  The inequality  $U_{(16)}(20) = 1.0000 > 0.9853 = v_{0.05}(5)$ (note that $U_{(16)}(20)$ corresponds to the fifth largest observation in absolute value) implies that  $d_1=5$. So it is possible that the number of outliers might be greater than 5. We reject the largest in absolute value 20th observation as an outlier and continue the search of outliers.

{\it Step 2.} The inequality $ U_{(15)}(19) = 1.0000 > 0.9853 = v_{0.05}(5)$ (note that $U_{(15)}(19)$ corresponds to the fifth largest observation in absolute value from the remaining 19 observations) implies that $d_2=5$. So it is possible that the number of outliers might be greater than 6. We declare the second largest in absolute value  observation as an outlier. So two observations (19th and 20th) are declared as outliers. We continue the search of outliers. 

 {\it Step 3.}  The inequality  $U_{(14)}(18) = 0.999997 > 0.9853 = v_{0.05}(5)$ implies that $d_3=5$.  We declare the third largest in absolute value  observation as an outlier. So three observations (2nd, 19th and 20th) are declared as outliers. We continue the search of outliers.

 {\it Step 4.} The inequalities $U_{(13)}(17) = 0.084290 < 0.9853 = v_{0.05}(5)$ and
$U_{(14)}(17) = 0.999940 > 0.9853 = v_{0.05}(5)$ imply that $d_4 = 4$. So four
additional observations (the fourth, fifth, sixth and seventh largest in absolute value observations), namely the 3d, 1st, 17th, and 7th are declared as outliers, The outlier search is finished. In all,   7 observations were declared as outliers: 1-3,17-20, as was expected. Note that  since the outlier search procedure was done after rejection of the null hypothesis, the significance level did not change.

We created R package \href{https://github.com/linas-p/outliersTests}{outliersTests}\footnote{
The R package outliersTest package can be accessed:
\href{https://github.com/linas-p/outliersTests}{https://github.com/linas-p/outliersTests}
} to be able to use the proposed BP test in practice within R package.

\section{Generalization of Davies-Gather outlier identification method }
\label{sec:davies}

Let us consider location-scale families.
 Following the idea of Davies-Gather \cite{davies_identification_1993}  define an empirical analogue
of the right outlier region  as a random region
\begin{equation}
OR_r(\alpha_n)=\{x:\,x>\hat\mu+\hat\sigma g_{n.\alpha}\},
\end{equation} where $g_{n.\alpha}$ is found using the condition
\begin{equation}\label{eq:g radimas}
\PP\{X_i\bar{\in}OR_r(\alpha_n),i=1,\ldots,n|H_0\}=1-\alpha,
\end{equation}
 and $\hat\mu,\hat\sigma$ are robust equivariant
estimators of the parameters $\mu,\sigma$.

Set $$
\hat Y_{(n)}=({X}_{(n)}-\hat\mu)/\hat\sigma.$$
The distribution of $\hat Y_{(n)}$ is parameter-free under $H_0$.

The equation \eqref{eq:g radimas} is equivalent to the equation
equation
$$
\PP\{\hat Y_{(n)}\leq g_{n,\alpha}\}|H_0\}=1-\alpha.
$$
So $g_{n,\alpha}$ is the upper 
$\alpha$ critical value of the random variable $\hat Y_{(n)}$. It is easily
computed  by simulation.
\vskip 0.2cm
\noindent {\it Generalized Davies-Gather method  for right outliers identification:} if $\hat Y_{(n)}\leq g_{n,\alpha}$, then it is concluded that right outliers are absent. The probability of such event is $\alpha$.
 If $\hat Y_{(n)}>g_{n,\alpha}$, then it is concluded that right outliers exist. 
 The value $x_i$ of the
random variable $X_i$ is admitted as an outlier if $x_i\in
OR_r(\alpha_n)$, i.e. if
$x_i>\hat\mu+\hat\sigma g_{n,\alpha}$. Otherwise it is admitted as a
non-outlier. \vskip 0.2cm

An empirical analogue
of the left outlier region  as a random region
\begin{equation}
OR_l(\alpha_n)=\{x:\,x<\hat\mu+\hat\sigma h_{n.1-\alpha}\},
\end{equation} where $h_{n.1-\alpha}$ is found using the condition
\begin{equation}\label{eq:g radimas 2}
\PP\{X_i\bar{\in}OR_l(\alpha_n),i=1,\ldots,n|H_0\}=1-\alpha,
\end{equation}
Set $$
\hat Y_{(1)}=({X}_{(1)}-\hat\mu)/\hat\sigma.$$
The distribution of $\hat Y_{(1)}$ is parameter-free under $H_0$.

The equation \eqref{eq:g radimas 2} is equivalent to the equation
equation
$$
\PP\{\hat Y_{(1)}\geq h_{n,1-\alpha}|H_0\}=1-\alpha.
$$
So $h_{n,\alpha}$ is the upper 
$1-\alpha$ critical value of the random variable $\hat Y_{(1)}$. It is easily
computed  by simulation.
\vskip 0.2cm
\noindent {\it Generalized Davies-Gather method  for left outliers identification:} if $\hat Y_{(1)}\geq h_{n,1-\alpha}$, then it is concluded that left outliers are absent. The probability of such event is $\alpha$.
 If $\hat Y_{(1)}<h_{n,\alpha}$, then it is concluded that left outliers exist. 
 The value $x_i$ of the
random variable $X_i$ is admitted as an outlier if $x_i\in
OR_l(\alpha_n)$, i.e. if
$x_i<\hat\mu+\hat\sigma h_{n,1-\alpha}$. Otherwise it is admitted as a
non-outlier. \vskip 0.2cm

Let us consider two-sided case. 

If the distribution of $X_i$ is symmetric, 
 then the empirical analogue
of the  outlier region  is the random region
\begin{equation}
OR(\alpha_n)=\{x:|x-\hat\mu|>\hat\sigma g_{n.\alpha/2}\}.
\end{equation}
In this case 
$$1-\alpha=\PP\{X_i\in OR(\alpha_n),i=1,...,n|H_0\}=\PP\{|\hat Y|_{(n)}\leq g_{n.\alpha/2}\}.
$$
{\it Generalized Davies-Gather method for left and right outliers identification (symmetric distributions):} if 
$|\hat Y|_{(n)}\leq g_{n.\alpha/2}$, then it is concluded that  outliers are absent. The probability of such event is $\alpha$.
 If $|\hat Y|_{(n)}> g_{n.\alpha/2}$, then it is concluded that  outliers exist. 
 The value $x_i$ of the
random variable $X_i$ is admitted as a left outlier if $x_i<\hat\mu-\hat\sigma g_{n,\alpha/2}$, it is admitted as a right outlier if
$x_i>\hat\mu+\hat\sigma g_{n,\alpha/2}$. Otherwise it is admitted as a
non-outlier.


If distribution of $X_i$ is non-symmetric, then   the empirical analogue
of the  outlier region  is defined as follows:
\begin{equation*}
OR(\alpha_n)=\{x\in {\bf R}/[\hat\mu+\hat\sigma g_{n,1-\alpha/2},\hat\mu+\hat\sigma g_{n,\alpha/2})]\},
\end{equation*}
In this case
$$
    1-\alpha=\PP\{X_i{\in}[\hat\mu+\hat\sigma h_{n,1-\alpha/2},\,\hat\mu+\hat\sigma g_{n,\alpha/2}],i=1,\ldots,n|H_0\}=
$$$$\PP\{h_{n,1-\alpha/2}\leq \hat Y_{(1)}\leq \hat Y_{(n)}\leq g_{n,\alpha/2})|H_0\}.
$$
\vskip 0.2cm
{\it Generalized Davies-Gather method for left and right outliers identification (non-symmetric distributions):} if $\hat Y_{(1)}\geq h_{n,1-\alpha/2}$
and $\hat Y_{(n)}\leq g_{n,\alpha/2}$, then it is concluded that  outliers are absent. The probability of such event is $\alpha$.
 If $\hat Y_{(1)}< h_{n,1-\alpha/2}$ or $\hat Y_{(n)}> g_{n,\alpha/2}$, then it is concluded that  outliers exist. 
 The value $x_i$ of the
random variable $X_i$ is admitted as a left outlier if $x_i<\hat\mu+\hat\sigma h_{n,1-\alpha/2}$, it is admitted as a right outlier if
$x_i>\hat\mu+\hat\sigma g_{n,\alpha/2}$. Otherwise it is admitted as a
non-outlier.

\section{Short survey of multiple outlier identification methods for normal data   }
\label{sec:compare}
\subsection{Rosner's method}  Let us formulate Rosner's method
in the form mostly used in practice. Suppose that the number of outliers does not exceed $s$ and
the two-sided alternative  is considered.  Set (see \cite{rosner_detection_1975,rosner1977percentage})
\begin{eqnarray*}
R_{1} =\max_{1\leqslant j \leqslant n} |\tilde Y_j| =
\max_{1\leqslant j \leqslant n} {|X_j - \bar{X}| }/{S_X}, \ \
\ S_X^2 = \sum_{j=1}^{n}(X_{(j)} - \bar{X})^2/{(n-1)}.
\end{eqnarray*}
$|\tilde Y_j|=|(X_j-\bar X)/S_X|$ may be interpreted as a distance
between $X_j$ and $\bar X$. Remove the observation $X_{j_1}$ which
is most distant from $\bar X$. This maximal  distance is $R_{1}$.
The value of $X_{j_1}$ is a possible candidate for contaminant.

Recompute the statistic using $n-1$ remaining observations and
denote by $R_{2}$ the obtained statistic. Remove the observation
$X_{j_2}$ which is most distant from the new empirical mean. The
value of $X_{j_2}$ is also possible candidate for contaminant.
Repeat the procedure until the statistics  $R_{1}, \cdots, R_{s}$
are computed. So we obtain all possible candidates for
contaminants. They are values of $X_{j_1},\ldots,X_{j_s}$

Fix $\alpha$ and find $\lambda_{in}$ such that 
\begin{equation*}
 \PP \{ R_{1} > \lambda_{in}  | H_0\}=...= \PP \{ R_{s} > \lambda_{in}  | H_0\},\quad \PP \{ \cup_{i=1}^{s} \{R_{i} > \lambda_{in} \} | H_0\} = \alpha.
\end{equation*}
If $n>25$, then the  approximations
\begin{eqnarray*}
\lambda_{in}\approx t_{\frac{\alpha}{2(
n-i-1)}}(n-i+1)\sqrt{\frac{n-i }{n-i-1+t^2_{\frac{\alpha}{2(
n-i-1)}}(n-i+1)}}\sqrt{1-\frac{1}{n-i+1}},
\end{eqnarray*}
are recommended  (see \cite{rosner_detection_1975}); here $t_{p}(\nu)$ is the $p$ critical value of the Student
distribution with $\nu$ degrees of freedom.
\vskip 0.2cm
 {\it Rosner's method for left and right outliers identification}: if $R_{i}\leqslant \lambda_{in}$ for all
$i=1,\cdots, s$, then it is concluded that outliers are
absent.
If there exists $i_0\in
\{1,\ldots,s\}$ such that $R_{i_0}> \lambda_{i_0n}$, i.e. the
event $ \cup_{i=1}^{s} \{R_{i} > \lambda_{in} \}$ occurs, then it is concluded that outliers exist. In this case 
classification of observations to outliers and non-outliers is done in the following way: if $R_{s}
> \lambda_{sn}$, then it is concluded that there are $s$ outliers and they are values of  $X_{j_1},\ldots,X_{j_s}$. If
$R_j \leqslant \lambda_{jn}$ for $j=s,s-1,\ldots,i+1$, and $R_{i} >
\lambda_{in}$,  then it is concluded that there are $i$
outliers and
 they are values of   $X_{j_1},\ldots,X_{j_i}$.

If   right outliers are searched, then define
 $R^+_{1} =\max_{1\leqslant i \leqslant n}\tilde Y_i,
$ and repeat the above procedure taking approximations
\begin{eqnarray*}
\lambda^+_{in} \approx t_{\frac{\alpha}{
n-i-1}}(n-i+1)\sqrt{\frac{n-i }{n-i-1+t^2_{\frac{\alpha}{
n-i-1}}(n-i+1)}}\sqrt{1-\frac{1}{n-i+1}}.
\end{eqnarray*}

\begin{figure}[t]
  \centering
  \includegraphics[width=0.7\linewidth]{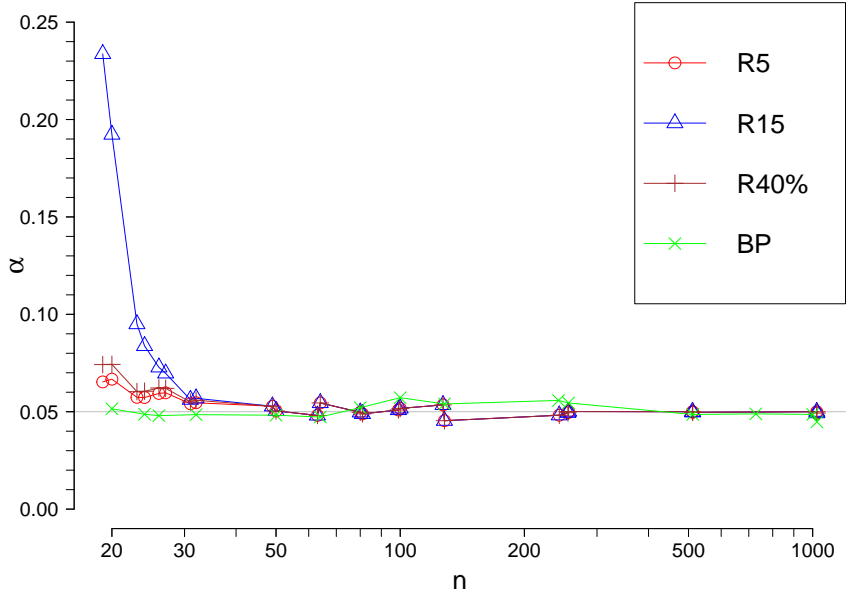}
    \caption{The true values of the significance level of Rosner's and BP tests in function of $n$ for different values of $s$ ($\alpha=0.05$ is used in approximations).  
	}\label{fig:hrejectionv2}

\end{figure}

Denote by $R_s$ the Rosner's test with a fixed upper limit  $s$. 
Our simulation results confirm  that  the true significance level is different from the level $\alpha$ suggested by the approximation  when $n$ is not large. Nevertheless, it is approaching
$\alpha$ as $n$ increases, see    Figure \ref{fig:hrejectionv2}. The true significance value of the BP test, which uses asymptotic values of the test statistic are 
also presented in Figure \ref{fig:hrejectionv2}.

\subsection{Bolshev's method}

Suppose that the number of contaminants does not
exceed $s$. For $i = 1,\cdots,n$ set
$$
\hat Y_i=(X_i-\bar X)/s,\quad \tau_i^+= n\cdot (1-T_{n-2}(\hat Y_i)) ,\quad \tau_i= n\cdot (1-T_{n-2}(|\hat Y_i|)),
$$
where $\bar X$ and $s$ are the empirical mean and standard deviation, $T_{n-2}(x)$ is the c.d.f. of Thompson's distribution with $n-2$ degrees of freedom. 

 Let us consider search for right outliers. Note that
 the
largest  $s$ observations $X_{(n-s+1)},...,X_{(n)}$ define the
smallest $s$ order statistics $\tau^+_{(1)}\leq...\leq
\tau^+_{(n)}$. Possible candidates for outliers are namely the
values of  $X_{(n-s+1)},...,X_{(n)}$.

Set $\tau^+ = \min_{1\leqslant i \leqslant s}{\tau_{(i)}^+}/{i}.$
\vskip 0.2cm \noindent {\it Bolshev's method for right outliers search}. If $\tau^+ \geq
\tau_{1-\alpha}^+(n,s)$, then it is concluded that outliers are
absent;
 here
$\tau_{1-\alpha}^+(n,s)$ is the $1-\alpha$ critical value of the
test statistic under $H_0$. If $\tau^+ < \tau_{1-\alpha}^+(n,s)$, then it is concluded that outliers exist. 
In such a case outliers  are
selected in the following way: if  ${\tau^+ _i}/{i} <
\tau_{1-\alpha}^+(n,s)$ then the value of the order statistic
$X_{(n-i+1)}$ is  admitted as an outlier, $i=1,...,s$. \vskip
0.2cm 

In the case of left and right outliers search  Bolshev's method uses
$\tau_{(i)}$ instead of $\tau^+_{(i)}$, defining the statistic
 $ \tau = \min_{1\leqslant i \leqslant
s}{\tau_{(i)}}/{i}.$
\vskip 0.2cm \noindent {\it Bolshev's method for left and right outliers search}. If $\tau \geq
\tau_{1-\alpha}(n,s)$, then it is concluded that outliers are
absent;
 here
$\tau_{1-\alpha}(n,s)$ is the $1-\alpha$ critical value of the
 statistic $\tau$ under $H_0$.  If $\tau < \tau_{1-\alpha}(n,s)$, then it is concluded that outliers 
exist. In such a case they are selected in
the following way: if  ${\tau _i}/{i} < \tau_{1-\alpha}(n,s)$
then the observation corresponding to $\tau_i$ is  admitted
as an outlier, $i=1,...,s$. \vskip 0.2cm

\subsection{Hawking's method }

Suppose that the number of contaminants does not
exceed $s$. Let us consider the search for right outliers.  For $k=1,...,s$ set
 \begin{equation*}
b^+_{k} = \frac{1}{\sqrt{k(n-k)}}\sum_{i=1}^k \tilde Y_{(n-i+1)}=
\frac{1}{\sqrt{k(n-k)}}\sum_{i=1}^k (X_{(n-i+1)}-
\bar{X})/S_X.
\end{equation*}
$b^+_{k}$ proportional to the sum of $k$ largest $\tilde Y_{(n-i+1)}$.  Set $B^+= \max_{1\leqslant k \leqslant s } b^+_{k}.$
\vskip 0.2cm \noindent {\it Hawking's method}. If $B^+ \leq
B^+_{\alpha}(n,s)$ then it is concluded that outliers are absent;
 here
$B^+_{\alpha}(n,s)$ is the $\alpha$ critical value of the 
statistic under $H_0$. If
$B^+ > B^+_{\alpha}(n,s)$, then it is concluded that outliers exist. In such a case  outliers are selected in the
following way: if  $b^+_{i} > B^+_{\alpha}(n,s)$, then the value of
the order statistic $X_{(n-i+1)}$ is  admitted as an outlier,
$i=1,...,s$.

\section{Comparative analysis of outlier identification methods by simulation}
\label{sec:simul}

In the case of location-scale classes  probability distribution of all considered test statistics does not depend on $\mu$ and $\sigma$, so we generated samples of various sizes $n$ with $n-r$
observations having the c.d.f. $F_0$  and  $r$ observations
having  various alternative distributions concentrated in the outlier region. We shall call such observations "contaminant outliers", shortly
$c-outliers$. As was mentioned, outliers which are not c-outliers, i.e. outliers from regular observations with the c.d.f. $F_0$,  are very rare.

 We repeated simulations $M=100000$ times and  using various methods we classified
observations to outliers and non-outliers and computed the mean number $ D_{O_cO}$ of correctly identified
c-outliers, the  mean number $ D_{ON}$ of  c-outliers which were not
 identified,   the mean number $ D_{NO}$ of  non c-outliers
admitted as outliers, and  the  mean number $ D_{NN}$ of non c-outliers admitted as
non-outliers.


An outlier identification method is ideal if each outlier is detected and each non-outlier is declared as a  non-outlier.
In practice it is impossible to do with the probability one. Two errors are possible: a) an outlier is not declared as such (masking effect); b)  a non-outlier is
declared as an outlier (swamping effect). We shall write shortly "masking value" for the mean number of non-detected c-outliers  and
"swamping value" for the mean number of "normal" observations declared as outliers in the simulated samples.

If  swamping is small for two tests  then  a test with smaller masking effect should be
preferred because in this case the distribution of the data remaining after excluding of suspected outliers should be closer
to the distribution of non-outlier data.

From the other side, if swamping for  Method 1 is considerably bigger
 than swamping of  Method 2 and  masking is smaller for Method 1, then  it does not
mean that Method 1 is better because this method rejects many extreme
non-outliers from the tails of the regular distribution $F_0$ and the sample remaining after
classification may be not treated as a sample from this regular
distribution even if all c-outliers are eliminated.

For various  families of distributions, sample sizes $n$, and alternatives we
 compared Davies-Gather (DG) and new (BP) methods performance. In the case of normal distribution
we also compared them with Rosner's, Bolshev's and Hawking's methods.

We used two different classes of alternatives:
in the first case c-outliers are spread widely in the outlier region around the mean, in the second case c-outliers are concentrated in a very
short interval laying in the outlier region. More precisely, if right outliers were searched, then 
we simulated $r$ observations concentrated in in the right
outlier region $out_r(\alpha_n,F_0)=\{x:x>x_{\alpha_n}\}$ using the following alternative families of distribution: 

1) Two parameter exponential distribution
$\mathcal{E}(\theta,x_{{\alpha_n}})$ with the scale parameter $\theta$.
 If $\theta$ is small, then outliers are concentrated near the border of the outlier region. If $\theta$  is large then  outliers are widely spread in the outlier region. If $\theta$ increases, then the mean of outlier distribution increases. Note that even if $\theta$ is very near 0 and the true number of outliers $r$ is large, these 
 outliers may corrupt strongly the data making tails of histogram two heavy. 

 2) Truncated normal distribution
$\mathcal{T}\mathcal{N}(x_{\alpha_n}, \mu, \rho)$ with the location and scale parameters $\mu,\rho$ ($\mu>x_{\alpha_n})$.  If $\rho$ is small then  this distribution is concentrated in a small interval around $\mu$. If $\mu$ increases, then the mean of outlier distribution increases. 

For lack of place we present a small part of our investigations. Note that the results are very
similar for all sample sizes $n\geq 20$. Multiple outlier problem is not very relevant for smaller sample sizes.

\subsection{Investigation of outlier identification methods  for normal data}

We  use notation  $B,H, R, DG$, and $BP$ for the Bolshev's, Hawking's, Rosner's, Davies-Gather's, and the new methods, respectively. If $DG$ method is based on maximum likelihood estimators, then we write $DG_{ml}$ method, if it is based on robust estimators, we write $DH_{rob}$ method.

For comparison of above considered methods we fixed the significance level $\alpha=0.05$.  We remind that the significance level $\alpha$ is the probability to reject minimum one observation as an outlier under the hypothesis $H_0$ which means that all observations are realizations of i.i.d. having the same normal distribution. The only test, namely R method uses approximate critical values of the test statistic, so the significance values for this test is only approximately $0.05$ and depends on $s$ and $n$. In Figure 1 the true significance level value for $s=5,15$ and $[0.4 n]$ in function of $n$ are given. 

The $B,H$,  and $R$ tests  methods have a drawback that the upper bound for the possible number of outliers $s$ must be fixed. The BP and DG tests have an advantage that they do not 
require it. 

Our investigations showed that H,B and $DG_{ml}$ methods have other serious drawbacks.  
 So firstly let us look closer at these methods.

\begin{figure}[t]
    \centering
  \includegraphics[width=0.6\linewidth]{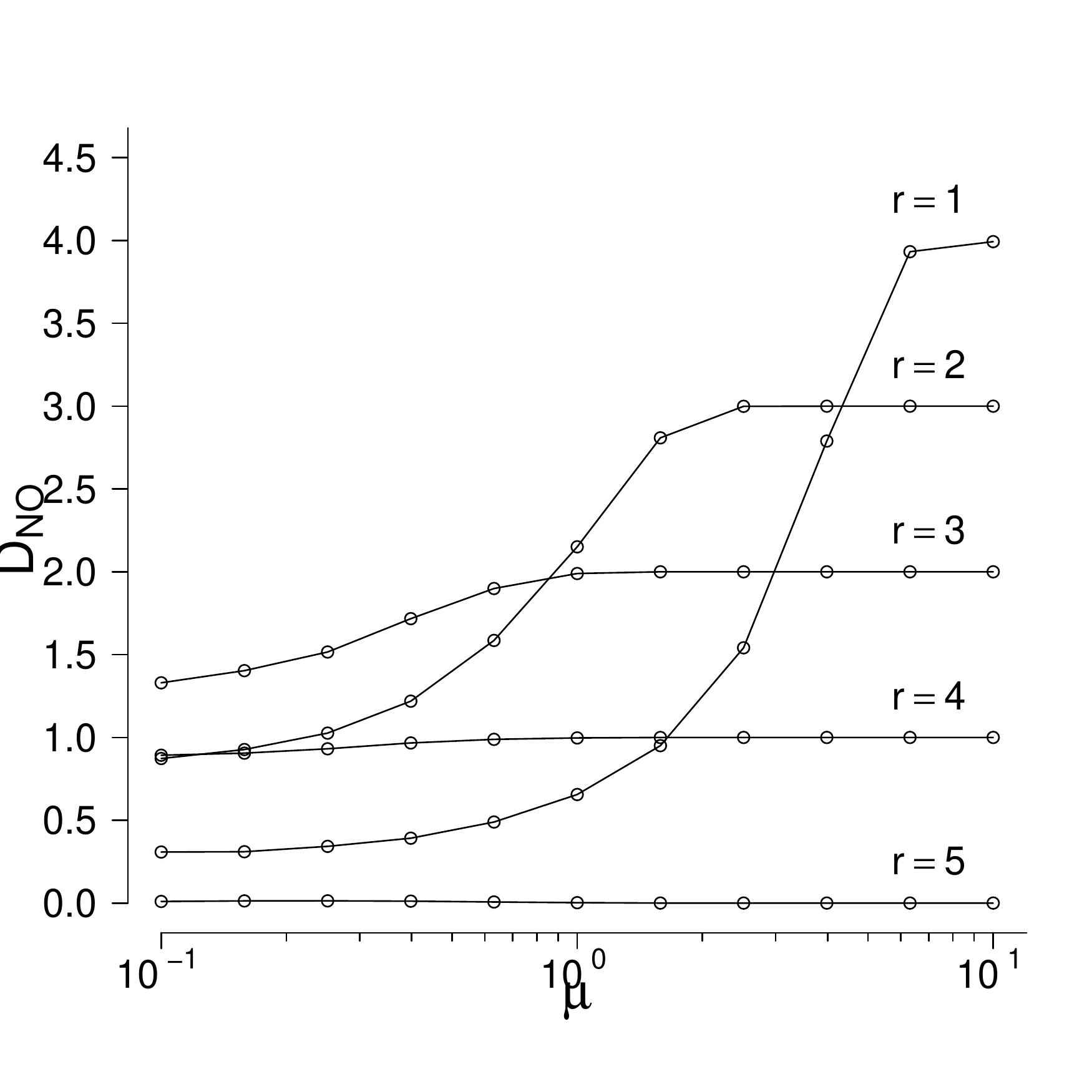}
  \caption{Hawkin's method: the values of $D_{NO} + D_{OO}$ in function of $\mu$ and $r$ ($n=100$, $s=5$).}\label{fig:hawkins1}

\end{figure}

\vskip 0.2cm

\begin{table}[ht]
\caption{Hawkin's method: the values of $D_{NO} + D_{OO}$ in function of $\mu$ and $r$ ($n=100$, $s=5$).}\label{tab:hawkins1}
\centering 
\begin{tabular}{lllll}
  \hline
r $\backslash$ $\mu$ & 0.1 & 1 & 6.3 & 10 \\ 
  \hline
 $1$ & 0.31 + 0.00 & 0.66 + 0.00 & 3.93 + 1.00 & 3.99 + 1.00 \\ 
$2$ & 0.87 + 0.00 & 2.15 + 0.06 & 3.00 + 1.21 & 3.00 + 2.00 \\ 
 $3$& 1.33 + 0.08 & 1.99 + 0.84 & 2.00 + 2.00 & 2.00 + 2.00 \\ 
$4$ & 0.89 + 0.58 & 1.00 + 1.42 & 1.00 + 3.00 & 1.00 + 3.00 \\ 
$5$ & 0.01 + 1.15 & 0.00 + 2.03 & 0.00 + 3.02 & 0.00 + 3.96 \\ 
   \hline
\end{tabular}

\end{table}

If the true number of c-outliers $r$ exceeds $s$, then the B and H methods can not find them even if they are very far from the limits of the outlier region. 
Nevertheless, suppose that $r$ does not exceed $s$ and look at the performance of the H method. Set $n=100$, $s=5$, and suppose that c-outliers are generated by 
right-truncated normal distribution $\mathcal{T}\mathcal{N}(x_{\alpha_n}, \mu, \rho)$ with fixed $\rho$ and increasing $\mu$. Note that the true number of c-outliers is supposed to be unknown but do not exceed $s=5$. In Figure \ref{fig:hawkins1} the mean numbers of rejected  non-c-outliers $D_{NO}$ are given in function of the parameter $\mu$  (the value of the parameter $\rho=0.1^2$ is fixed) for fixed values of  $r$ see Figure \ref{fig:hawkins1}. In Table \ref{tab:hawkins1}  the values of $D_{NO}$ plus the values of the mean numbers of truly rejected c-outliers are given.   The Table \ref{tab:hawkins1} shows that if  $r=1$, then if $\mu$ is sufficiently large, the c-outlier is found but the number of rejected non-c-outliers $D_{NO}$ increases to 4, so swamping is very large. Similarly, if $r=2$, then  $D_{NO}$ increases to 3, so swamping is large. Beginning from $r=3$ not all c-outliers are found even for large $\mu$.  Swamping is smallest if the true value $r$  coincides with $s$ but even in this case one c-outlier is not found even for large $\mu$. 
 Taking into account that the true number $r$ of c-outliers is not known in real data, the performance of the H methos is very poor. Results are similar for other  values of $n$, $s$, and  distributions of c-outliers. As a rule, H mehod finds rather well the c-outliers but  swamping is very large because this method has a tendency to reject a number near $s$ of observations for remote alternatives.  which is good if $r=s$ but is bad if $r$ is different from $s$.  

The B and $DG_{ml}$ tests have a drawback that they use maximum likelihood estimators which are not robust and estimate parameters badly in presence of outliers.  
Once more, set $n=100$, $s=5$, and suppose that c-outliers are generated by 
two-parameters exponential distribution $\mathcal{T}\mathcal{E}(x_{\alpha_n}, \theta)$ with increasing $\theta$. Swamping values are negligible in, so only masking values( mean numbers of non-rejected  c-outliers $D_{ON}$) are important. In  Figure  3 the masking values  in function of the parameter $\theta$  are given  for fixed values of  $r$. 

Both methods perform very similarly. The  masking values are   large for every value of $r>1$. If $r$ increases, then masking values increase, too. For example, if $r=5$, then almost 3 c-outliers from 5  are not rejected in average even for large values of $\theta$. 

Similar results hold taking other values of $n$, $s$ and various distributions of c-outliers.

The above analysis shows that the  B, H, $DG_{ml}$ methods have serious drawbacks,  so we exclude these methods from further consideration. 

\begin{figure}[t]

\centering
  \includegraphics[width=0.5\linewidth]{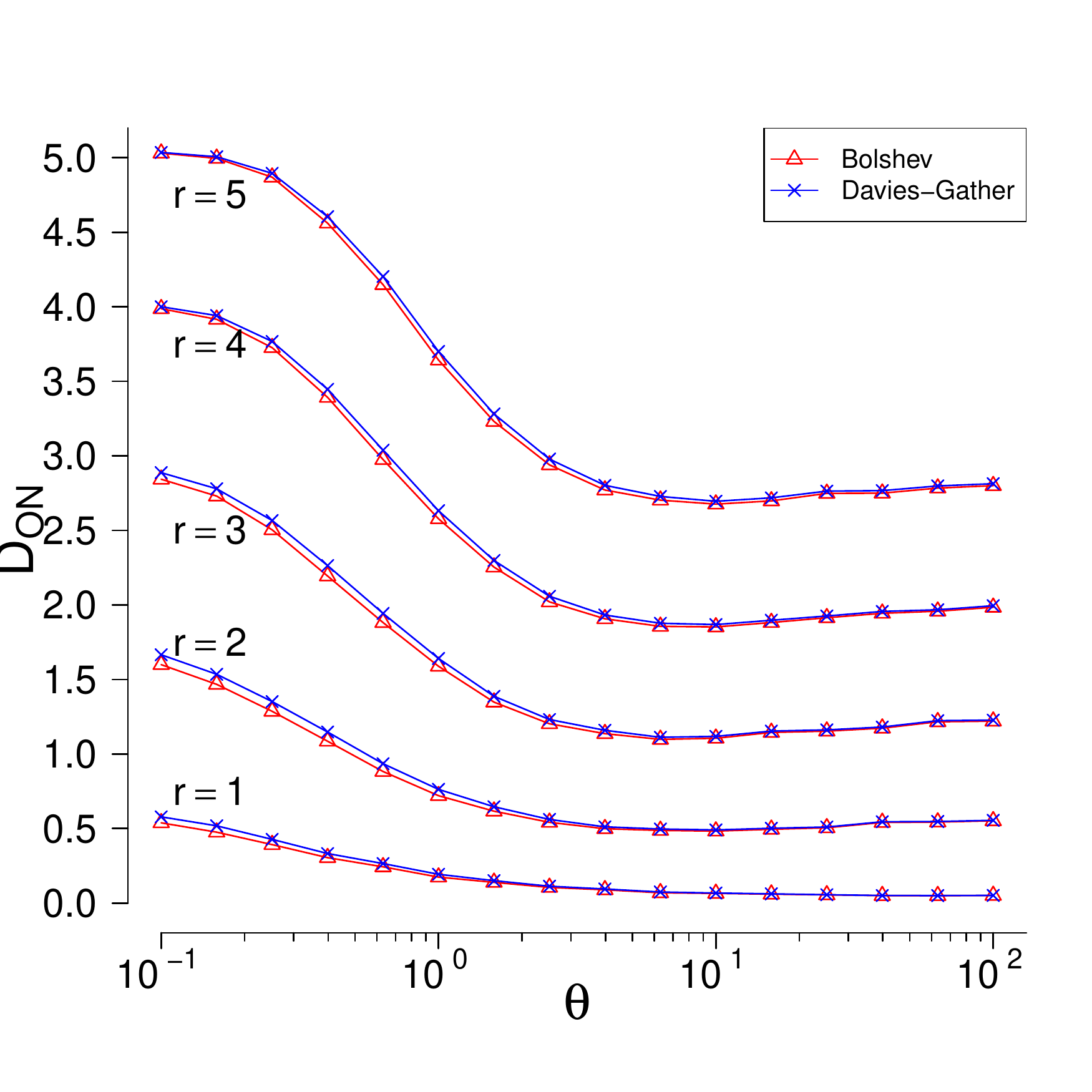}
	\caption{The number of outliers rejected as non-outliers ($D_{ON}$). The alternative: two-sided, the outliers generated by two-parameters exponential distribution on both sides. }\label{fig:dgb}
\end{figure}

\vskip 0.2cm

Let us consider  the remaining three methods: R, DG, and BP. For small $n$ the true significance level of Rosner's test differ considerably from the
 suggested, so we present comparisons of tests performance for $n=50,100,1000$ (see Tables \ref{tab:norm50}-\ref{tab:norm1000}). Truncated exponential distribution was used for outliers simulation. Remoteness
 of the mean of outliers from the border of the outlier region is characterized by the parameter $\theta$.

  \begin{table}[h]
	\caption{The masking values $D_{ON}$ ($n=50$ and $n=100$).}\label{tab:norm50}
	\centering
\begin{tabular}{|l|l|lllll|}
  \hline
  & $n = 50$ & &&&& \\
  \hline
r&Method$\backslash\theta$ & 0.1 & 0.4 & 1 & 4 & 10 \\
  \hline
2&Rosner\_5 & 1.36 & 0.95 & 0.51 & 0.15& 0.06 \\
&  Rosner\_15 & 1.36 & 0.95 & 0.51 & 0.15& 0.06 \\
&  Rosner\_[0.4n] & 1.36 & 0.95 & 0.51 & 0.15& 0.06 \\
&  DGr\_rob & 1.56 & 1.17 & 0.71 & 0.24 & 0.10 \\
&  BP & 0.92 & 0.66 & 0.37& 0.10 & 0.04 \\
   \hline
5&Rosner\_5 & 3.79 & 3.31 & 2.11 & 0.48 & 0.16 \\
&  Rosner\_15 & 3.66 & 3.21 &2.04 & 0.46 & 0.16 \\
&  Rosner\_[0.4n] & 3.66 & 3.21 &2.04 & 0.46 & 0.16 \\
&  DG\_rob &4.70 & 4.10& 2.90 & 1.09 & 0.48 \\
&  BP & 2.00 & 1.68 & 1.18 & 0.40 & 0.15 \\
   \hline
8&Rosner\_5 & 8.00& 7.97 & 7.54 & 3.70 & 3.06 \\
&  Rosner\_15 &  5.70& 5.48 & 4.52 & 1.00 & 0.29 \\
&  Rosner\_[0.4n]& 5.70& 5.48 & 4.52 & 1.00 & 0.29 \\
&  DG\_rob & 7.90 & 7.49 & 6.10 & 2.67 & 1.24 \\
&  BP & 4.27 & 3.84 & 3.25 & 1.47 & 0.57 \\
   \hline
\end{tabular}\begin{tabular}{|l|lllll|}
  \hline
    & $n = 100$ & &&& \\
  \hline 
r & 0.1 & 0.4 & 1 & 4 & 10 \\
  \hline
2 & 1.19 & 0.71 & 0.33 & 0.09 & 0.04 \\
 &  1.19 & 0.71 & 0.33 & 0.09 & 0.04 \\
 &  1.19 & 0.71 & 0.33 & 0.09 & 0.04 \\
& 1.31 & 0.84 & 0.44 & 0.13 & 0.06 \\
 & 0.50 & 0.32 & 0.15 & 0.04 & 0.02 \\
   \hline
5 & 3.52 & 2.57 & 1.27 & 0.27 & 0.10 \\
 & 3.43 & 2.52 & 1.24 & 0.26 & 0.10 \\
 & 3.43 & 2.52 & 1.24 & 0.26 & 0.10 \\
 & 4.23 & 3.01 & 1.81 & 0.57 & 0.25 \\
 & 0.78 & 0.60 & 0.43 & 0.15 & 0.07 \\
   \hline
10 & 10.0  & 9.90 & 8.21 & 5.10 & 5.00 \\
 & 6.88 & 6.54 & 4.36 & 0.69 & 0.22 \\
 & 6.88 & 6.54 & 4.36 & 0.69 & 0.22 \\
 & 9.74 & 8.38 & 5.78 & 2.12 & 0.92 \\
 & 2.21 & 1.90 & 1.73 & 0.74 & 0.30 \\
   \hline
\end{tabular}

\end{table}

\begin{table}[h]
	\caption{The masking values $D_{ON}$ ($n=1000$).}\label{tab:norm1000}
\centering
\begin{tabular}{|l|l|lllll|}
  \hline
r&Method$\backslash\theta$ & 0.1 & 0.4 & 1 & 4 & 1000 \\
  \hline
5&Rosner\_5 & 2.15 & 0.69 & 0.29 & 0.07 & 0.00 \\
 & Rosner\_15 & 2.12 & 0.66 & 0.27 & 0.07 & 0.00 \\
 & Rosner\_[0.4n] & 2.12 & 0.66 & 0.27 & 0.07 & 0.00 \\
 & DG\_rob & 1.99 & 0.78 & 0.35
 & 0.09 & 0.00 \\
 & BP &0.25 & 0.23 & 0.22 & 0.11 & 0.00 \\
   \hline
20&Rosner\_5 & 19.0 & 15.8 & 15.0 & 15.0 & 15.0 \\
 & Rosner\_15 & 19.2 & 10.9 & 5.52 & 5.00 & 5.00 \\
 & Rosner\_[0.4n] & 12.7 & 6.94 & 1.76 & 0.30 & 0.00 \\
&  DG\_rob & 14.8 & 6.97 & 3.32 & 1.93 & 0.00 \\
&  BP & 0.29 & 0.26 & 0.23 & 0.18  & 0.00 \\
   \hline
100&Rosner\_5 & 100 & 99.9 & 96.7 & 95.0 & 95.0 \\
 & Rosner\_15 & 100  & 99.92 &  96.4  & 85.0 & 85.0\\
 & Rosner\_[0.4n]  & 55.8 & 56.8 & 50.4 & 4.43 & 0.01 \\
 & DG
\_rob & 100 & 89.9 & 61.6 & 22.2 & 0.1 \\
 & BP & 4.72 & 4.00 & 3.95 &3.58 & 0.04 \\
   \hline
\end{tabular}

\end{table}

 Swamping values $D_{NO}$  (the mean numbers of non-c-outliers declared as outliers) are very small for all tests. For example, even if $n=1000$, the  R and DG methods reject in average  as outliers only  0.05 from $n-r=995, 980, 900$  non-c-outliers. For the BP method this number is $0.25, 0.19, 0.05$ from $995, 980$, and 900 non-c-outliers, respectively. So only masking values $D_{ON}$ (the mean numbers of c-outliers declared as non-outliers) are important for outlier identification methods comparison.

 Necessity to guess the upper limit $s$ for a possible number of outliers is considered as  a drawback of the Rosner's method. Indeed, if the
 true number of outliers $r$ is greater than the chosen upper limit $s$, then $r-s$ outliers are not identified with the probability one.  And
  even if $r\leq s$, it is  not clear how important is closeness of $r$ to $s$. So first we investigated the problem of the upper limit choice.

  Here we present masking values $D_{ON}$ of the  Rosner's tests for $s=5,15$ and $[0.4n]$. Similar results are obtained for other values of $s$.

  Our investigations show that it is sufficient to fix $s=[0.4n]$, which  is clearly larger than it can be expected in real data. Indeed, Tables \ref{tab:norm50}-\ref{tab:norm1000} show that for $r>s$ $Rosner_5$ and $Rosner_{15}$   do not find $r-s$ outliers
  even if they are very remote, as it should be. Nevertheless, we see that even if the true number of outliers $r$  is much smaller than $[0.4n]$,  for any  considered $n$, $r\leq s=5,15$ the masking values of the $Rosner_{[0.4n]}$ test are  approximately the same (even a little smaller) as  the masking values of  the tests $Rosner_5$ and $Rosner_{15}$,  for $r>s$  they  are clearly smaller.

  Hence, $s=[0.4n]$ should be recommended for Rosner's test application, and performance of $Rosner_{[0.4]}$, Davies-Gather robust ($DG_{rob})$ and the proposed $BP$ methods should be compared.

  All three methods find all c-outliers if they are sufficiently remote.
  For $n=50$ the BP method gives uniformly smallest masking values and the DG method gives uniformly largest masking values for any considered $r$ in all diapason of alternatives. For $n=100$ and $r=2,5$ the result is the same. For
  $n=100$ and $r=10$ (it means that even for very small $\theta$ the data is seriously corrupted) the BP method is also the best except that for the  most remote alternatives the $Rosner_{[0.4n]}$ method slightly outperforms the BP method. For $n=1000$ and the most of alternatives the BP method strongly outperforms other methods, except the most remote alternatives.  
  
  The DG and Rosner's methods have very large masking if many outliers are concentrated near the outlier region border. In this case data is seriously corrupted, however,  
  these methods do not see outliers. 

Conclusion: in most considered situations  the BP method is the best outlier identification method. The second is Rosner's method with $s=[0.4]$, and the third is the Davies-Gather method based on robust estimation. Other methods have poor performance.

\subsection{Investigation of outlier identification methods  for other location-scale models}

\begin{figure}[t]
\centering
\subfloat[Logistic distribution.]{%
\resizebox*{5cm}{!}{\includegraphics{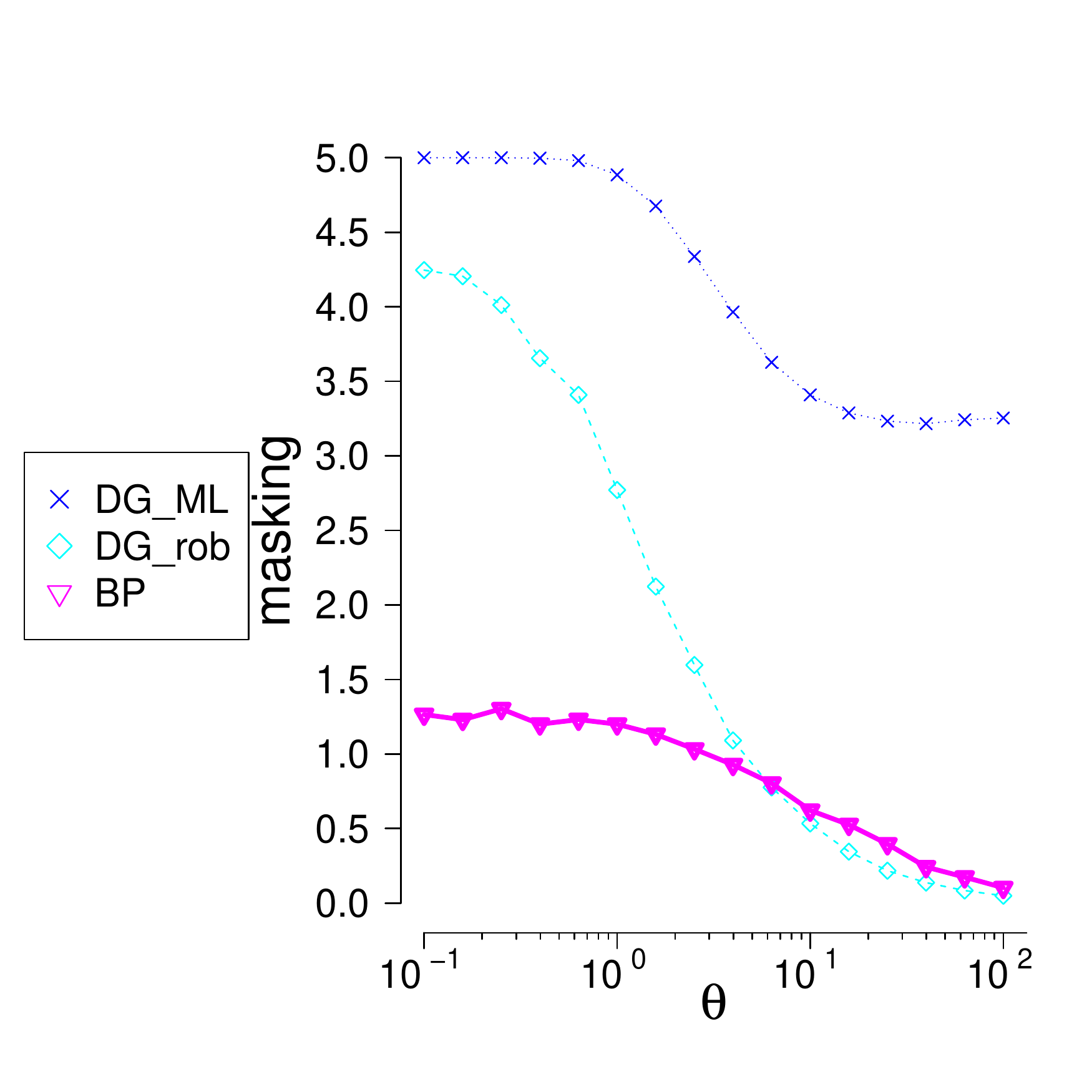}}}\hspace{5pt}
\subfloat[Laplace distribution.]{%
\resizebox*{5cm}{!}{\includegraphics{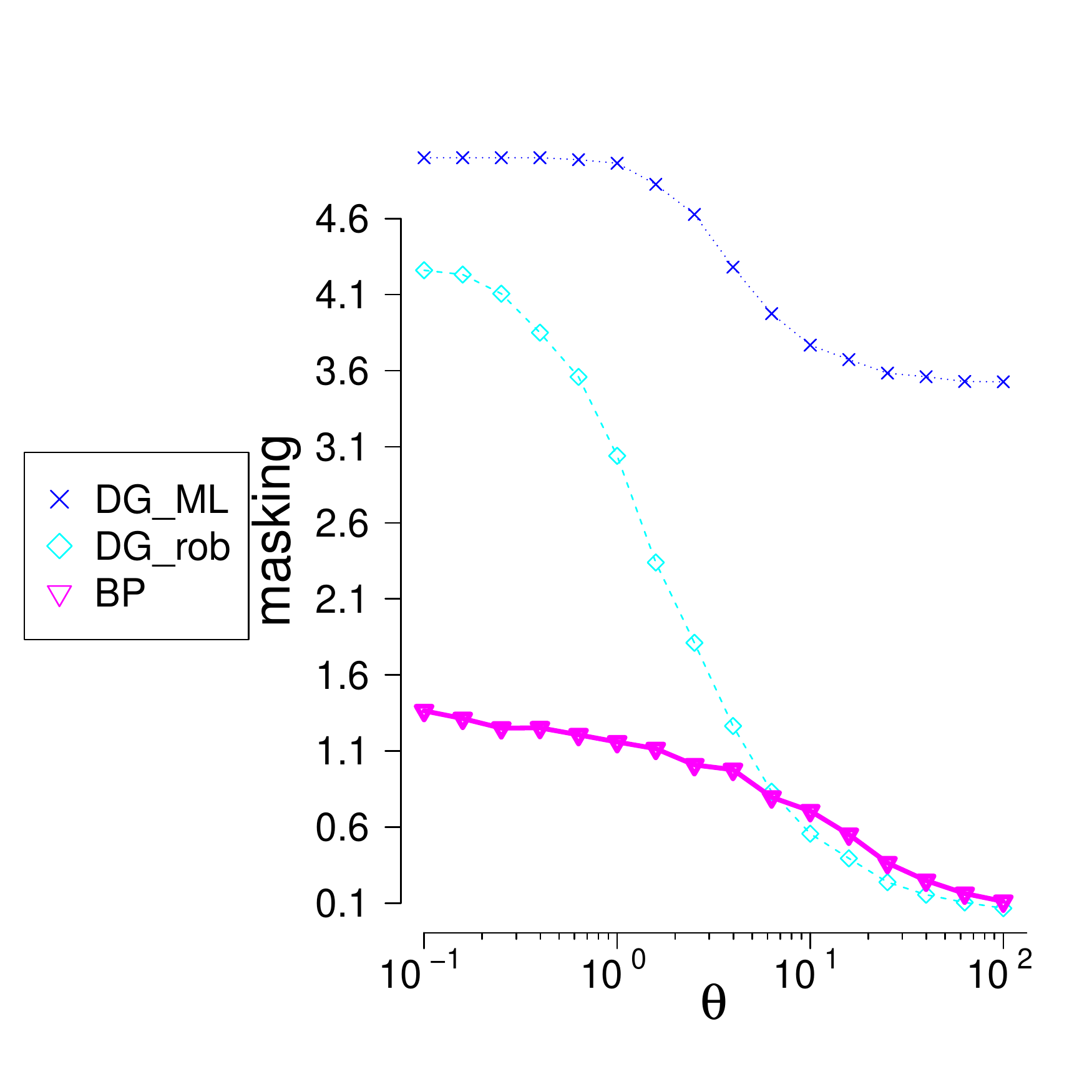}}}

	\caption{The difference between number outliers and rejected observations   given that sample size $n=100$ and $r=10$ outliers. }\label{fig:logist}
\end{figure}

\begin{figure}[t]
\centering
\subfloat[Extreme value II distribution.]{%
\resizebox*{5cm}{!}{\includegraphics{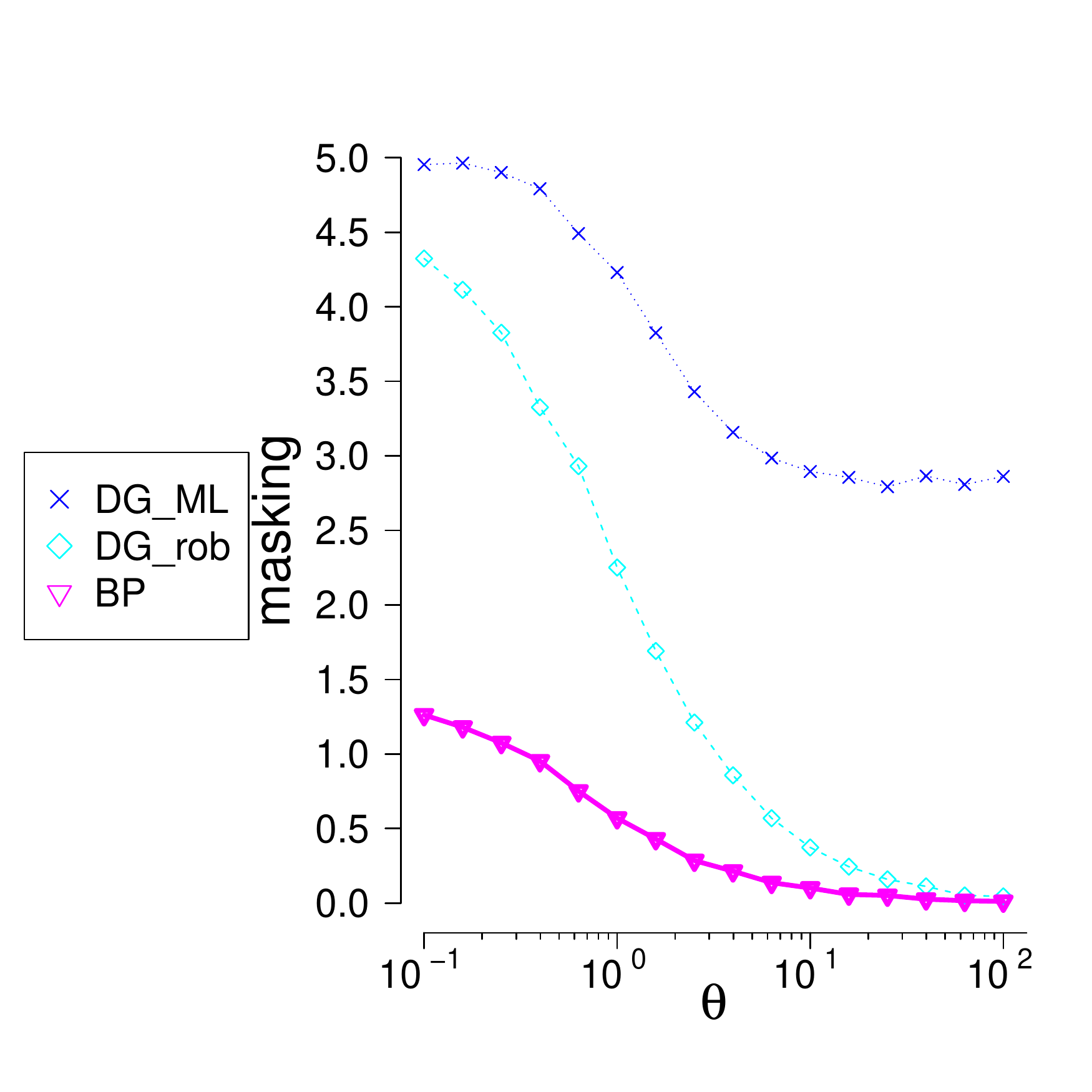}}}\hspace{5pt}
\subfloat[Cauchy distribution.]{%
\resizebox*{5cm}{!}{\includegraphics{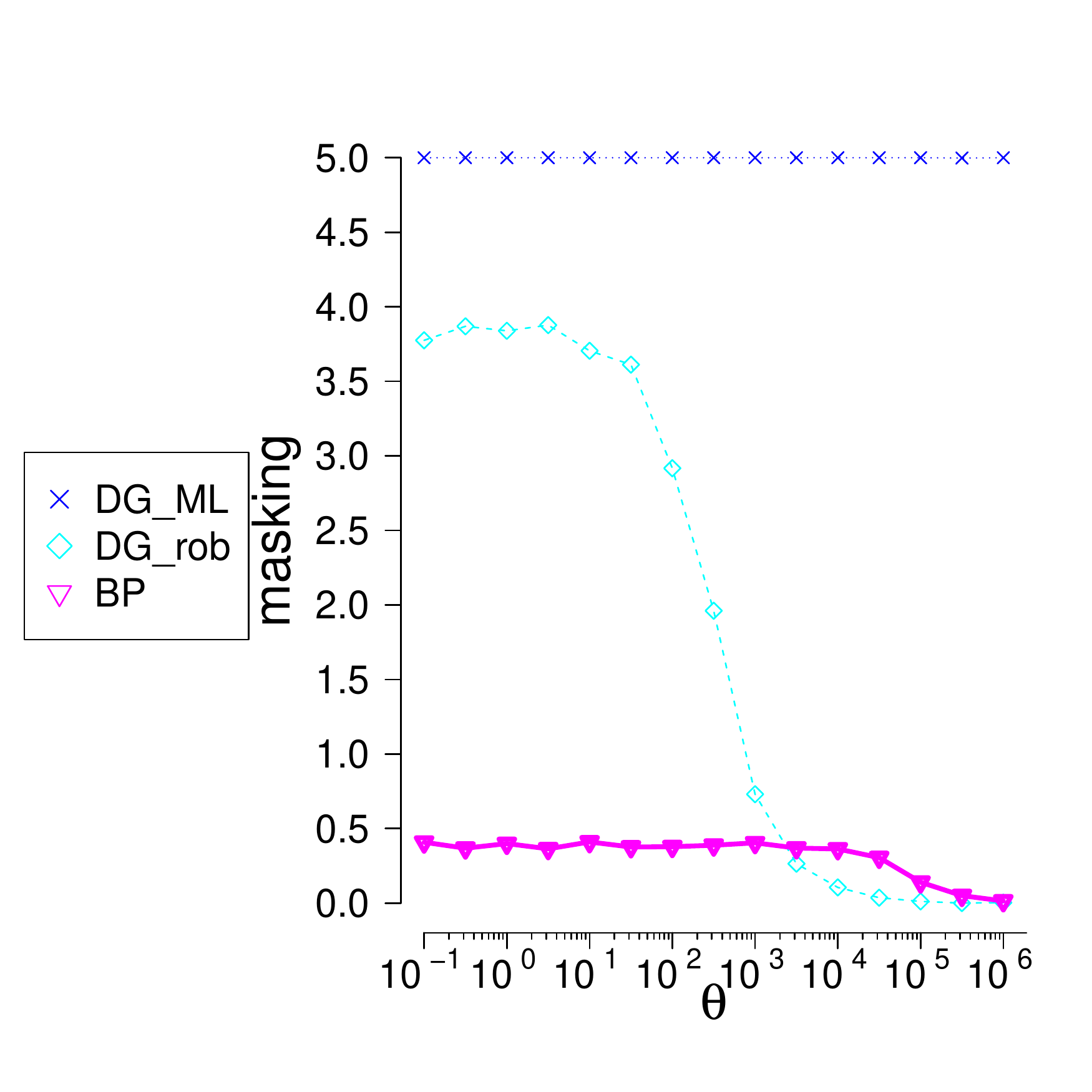}}}

	\caption{The difference between number outliers and rejected observations   given that sample size $n=100$ and $r=10$ outliers. }\label{fig:cauchy}
\end{figure}

We investigated performance of the new method for location-scale families different from normal. We compare 
the BP method with the generalized Davies-Gather method for logistic, Laplace (symmetric, $F\in {\cal G}_0\cap {\cal F}_{ls}$), extreme values 
(non-symmetric $F\in {\cal G}_0\cap {\cal F}_{ls}$), and Cauchy (symmetric, $F\in {\cal G}_1\cap {\cal F}_{ls}$) families. C-outliers were generating using truncated 
exponential distribution concentrated in two-sided outlier region.  Swamping values being small, masking value, see Table \ref{tab:other} and differences between the true number of c-outliers and the number of rejected observations, see Figures 
\ref{fig:logist}-\ref{fig:cauchy}, were compared. The BP and $DG_{rob}$ methods find very well the most remote outliers, meanwhile, the BP method identifies much better  closer outliers.  The $DG_{rob}$ method identifies badly multiple outliers concentrated near the border of the outlier region, whereas the BP method does well. The $DG_{ML}$ is not appropriate for multiple outlier search.

\begin{table}[ht]
\caption{Masking values for logistic, Laplace, extreme value II and Cauchy distribution, when $n=100$, $r=5$.}\label{tab:other}
\centering 
\begin{tabular}{|l|clll|clll|}
  \hline
  & Logistic & & & & Laplace & & &\\ \hline 
Method$\backslash \theta$ & 0.1 & 1 & 6.3 & 10 & 0.1 & 1 & 6.3 & 10 \\  
  \hline
DaviesGather\_ML & 5 & 4.89 & 3.64 & 3.42 & 5 & 4.96 & 3.98 & 3.78 \\ 
  DaviesGather\_Robust & 4.21 & 2.69 & 0.76 & 0.51 & 4.27 & 2.98 & 0.87 & 0.59 \\ 
  BP & 1.3 & 1.13 & 0.78 & 0.64 & 1.31 & 1.21 & 0.8 & 0.66 \\
   \hline
  \hline
   & Extreme value II &&&& Cauchy &&&\\ \hline 
Method$\backslash \theta$ & 0.1 & 1 & 6.3 & 10 & 1 & 100 & 1000 & $10^5$ \\ 
  \hline
DaviesGather\_ML & 4.96 & 4.19 & 3 & 2.9 & 5 & 5 & 5 & 5 \\ 
  DaviesGather\_Robust & 4.29 & 2.25 & 0.59 & 0.4  & 3.81 & 2.89 & 0.8 & 0.01 \\ 
  BP & 1.25 & 0.56 & 0.14 & 0.11 & 0.38 & 0.4 & 0.39 & 0.13 \\ 
   \hline
\end{tabular}

\end{table}

\newpage

\section{Conclusion}
In many situations, the proposed  outlier identification method  has superior performance as compared to  existing methods. The method is based on an asymptotic result, so it should be not applied for samples of very small size $n\leq 15$.

%
%
%
%
%
The formulated approach could be extended and used for  probability 
distribution families different from location-scale and shape-scale.

The R package outliersTests was created for the practical usage of proposed test.

\section*{Disclosure statement}
None

\section*{Funding}
This research did not receive any specific grant from funding agencies in the public, commercial, or non-profit sectors.

\bibliographystyle{plain}
\bibliography{Outliers}

\appendix
\section{Proofs of the theorems }

\noindent{\bf Proof of Theorem 3.1.}
Note that
$$
\frac{\hat Y_{(n-i+1)}-b_n}{a_n}=\frac{Y_{(n-i+1)}-b_n}{a_n}\,
\frac{\sigma}{\hat\sigma}+\frac{(\mu-\hat\mu)}{\hat\sigma a_n}+\frac{b_n}{a_n}\frac{\sigma-\hat\sigma}{\hat\sigma}
$$
The $s$-dimensional random vector such that its  $i$th component is the first  term of the right side  converges in distribution to the random vector given in the formulation of the theorem. It follows from Theorem 2.1.1 of \cite{de2007extreme} and \hyperref[cond:cond]{Condition A} \ref{cond:a}. So it is sufficient to show that  the second and the third terms converge to zero in probability. The second term is
$$
-\sqrt{n}f_0(F_0^{-1}(1-\frac1n))\frac{\sqrt{n}(\hat\mu-\mu)}{\hat\sigma},
$$
the third term is $$
-\sqrt{n}F_0^{-1}(1-\frac1n)f_0(F_0^{-1}(1-\frac1n))\frac{\sqrt{n}(\hat\sigma-\sigma)}{\hat\sigma}.
$$
By \hyperref[cond:cond]{Condition A} \ref{cond:c}
$$
\lim_{n\to \infty}\sqrt{n}F_0^{-1}(1-\frac1n)f_0(F_0^{-1}(1-\frac1n)))=\lim_{x\to \infty}\frac{xf_0(x)}{\sqrt{1-F_0(x)}}= 0.
$$
It also implies $\lim_{n\to\infty}\sqrt{n}f_0(F_0^{-1}(1-\frac1n))=0$  because $\lim_{n\to\infty}F_0^{-1}(1-\frac1n)=\infty.$  These results and \hyperref[cond:cond]{Conditions A}  \ref{cond:a}, \ref{cond:b} 
  imply the statement of the theorem.

\vskip 0.2cm
\noindent{\bf Proof of Theorem 3.2.}

For any $i=1,...,s$ the following equality holds:
\begin{equation}\label{pagrindine}
     \frac{|\hat Y|_{(n-i+1)}-b_{2n}}{{a_{2n}}}=\frac{|\hat Y|_{(n-i+1)}-|Y|_{(n-i+1)}}{{a_{2n}}}+\frac{|Y|_{(n-i+1)}-b_{2n}}{{a_{2n}}}.
\end{equation}
 The c.d.f. of  the random
variables $|Y_i|$ is $2F_0(x)-1$, so if $F_0\in {\cal G}_\gamma$, $\gamma\geq 0$ then $2F_0-1\in {\it \cal G}_\gamma$, and for the sequence $|Y_n|$ the normalizing sequences are $a_{2n},b_{2n}$.  So 
the $s$-dimensional random vector such that its  $i$th component is the second  term of the right side  converges in distribution to the random vector given in the formulation of the theorem. It follows from Theorem 2.1.1 of \cite{de2007extreme}. So it is sufficient to show that  the first term converges to zero in probability.
 
 Note that 
$|\hat Y_i|\leq |Y_i|+|\hat Y_i-Y_i|$,
 and 
$$
|\hat Y_i-Y_i|=
\frac1{\hat\sigma}|\mu-\hat\mu+(\sigma-\hat\sigma)Y_i|\leq \frac{|\hat\mu-\mu|}{\hat\sigma}+\frac{|\sqrt{n}(\hat\sigma-\sigma)|}{\hat\sigma}\frac1{\sqrt{n}}|Y|_{(n)}.
$$
So
$|\hat Y|_{(n-j+1)}\leq |Y|_{(n-j+1)}+\frac{|\hat\mu-\mu|}{\hat\sigma}+\frac{|\sqrt{n}(\hat\sigma-\sigma)|}{\hat\sigma}\frac1{\sqrt{n}}|Y|_{(n)}.$
Analogously, the inequality $ Y_i|\leq ||\hat Y_i|+|\hat Y_i-Y_i|$ implies that $|Y|_{(n-j+1)}\leq |\hat Y|_{(n-j+1)}+\frac{|\hat\mu-\mu|}{\hat\sigma}+\frac{|\sqrt{n}(\hat\sigma-\sigma)|}{\hat\sigma}\frac1{\sqrt{n}}|Y|_{(n)}.$

Theorem 2.1.1 in \cite{de2007extreme} applied to the random variables $|Y_i|$ implies that there exist 
a random variable $V_1$ with the c.d.f. $G(x)=e^{-e^{-x}}$ ($\gamma=0$) or  $G(x)=e^{-(-x)^{-1/\gamma}}, x<0$, $G(x)=1$, $x\geq 0$ ($\gamma>0$), such that 
\begin{equation} \label{eq:epsilonai}
\frac1{\sqrt{n}}|Y|_{(n)}=(b_{2n}+a_{2n}(V_1+o_P(1)))/\sqrt{n}.
\end{equation}
$$
|\frac{|\hat Y|_{(n-i+1)}-|Y|_{(n-i+1)}}{{a_{2n}}}|\leq \frac{|\sqrt{n}(\hat\mu-\mu)|}{\hat\sigma\sqrt{n}a_{2n}}+\frac{|\sqrt{n}(\hat\sigma-\sigma)|}{\hat\sigma}(\frac{b_{2n}}{\sqrt{n}a_{2n}}+\frac{V_1+o_p(1)}{\sqrt{n}}).
$$
The convergence $b_n\to\infty$ and  \hyperref[cond:cond]{Condition A} \ref{cond:c} imply:
$$
\lim_{n\to\infty}\frac{b_{2n}}{\sqrt{n}a_{2n}}=\lim_{n\to\infty}\sqrt{n}F_0^{-1}(1-\frac1{2n})f_0F_0^{-1}(1-\frac1{2n}))=$$
$$\frac1{\sqrt{2}}\lim_{x\to\infty}\frac{xf_0(x)}{\sqrt{1-F_0(x)}}=0,\quad \lim_{n\to\infty}\frac{1}{\sqrt{n}a_{2n}}=0.
$$
These results and \hyperref[cond:cond]{Conditions A} \ref{cond:a}, \ref{cond:b}   imply that the first term at the right of  (\refeq{pagrindine}) converges to zero in probability. 

\vskip 0.2cm
\noindent{\bf Proof of Theorem 3.3.}
The result follows from the equality
$$
\frac{\hat Y_{(n-i+1)}- b_n}{ a_n}=\frac{Y_{(n-i+1)}-b_n}{a_n}\frac{\sigma}{\hat\sigma}+\frac{b_n}{a_n}(\frac{\sigma}{\hat\sigma}-1)+\frac1{a_n}\frac{\mu-\hat \mu}{\hat\sigma}
,
$$
equivariance of the estimators $\hat\mu$, $\hat\sigma$ and the fact that the distribution of the random vector $(Y_1,...,Y_n)^T$ is parameter-free. 

\end{document}